\documentclass[12pt]{article}
\usepackage{amssymb,amsmath,amsthm,amsfonts,enumerate, bbm, mathdots}

\usepackage{latexsym}
\usepackage{amscd}
\usepackage{stmaryrd}
\usepackage{color}
\usepackage[all]{xypic}
\usepackage{epsfig}
\usepackage{graphics}
\usepackage{ifthen}
\usepackage{varioref}
\usepackage{rotating}
\usepackage{extarrows}
\usepackage{cite}
\usepackage{mathrsfs}

\numberwithin{equation}{section}

\theoremstyle{plain}
\newtheorem{thm}{Theorem}[section]
\newtheorem{cor}[thm]{Corollary}
\newtheorem{lem}[thm]{Lemma}
\newtheorem{prop}[thm]{Proposition}
\newtheorem{rem}[thm]{Remark}
\newtheorem{Def}[thm]{Definition}
\newtheorem{eg}[thm]{Example}

\definecolor{darkgreen}{rgb}{0.0625,0.64,0.0625}
\usepackage[%dvipdfm,
            pdfstartview=FitH,
            CJKbookmarks=true,
            bookmarksnumbered=true,
            bookmarksopen=true,
            colorlinks,
            pdfborder=001,
            linkcolor=blue,
            anchorcolor=green,
            citecolor=red
            ]{hyperref}

\newfont{\scyr}{wncyr10 scaled 550}

\def\proof{\noindent {\bf Proof.\;}}

\allowdisplaybreaks

\begin{document}

\title{Rota-Baxter systems and skew trusses}

\date{\small ~ \qquad\qquad School of Mathematical Sciences, Tongji University \newline No. 1239 Siping Road,
Shanghai 200092, China}

\author{Zhonghua Li\thanks{E-mail address: zhonghua\_li@tongji.edu.cn} ~and ~Shukun Wang\thanks{E-mail address: 2010165@tongji.edu.cn}}

\maketitle

\begin{abstract}
 As a generalization of skew braces, the notion of skew trusses was introduced by T. Brzezi\'{n}ski. It was shown that every Rota-Baxter group has the structure of skew braces by V. G. Bardakov and V. Gubarev. To investigate an analogue of Rota-Baxter groups which has the structure of skew trusses, we define Rota-Baxter systems. We study the relationship between Rota-Baxter systems and Rota-Baxter groups. Furthermore, we prove that a Rota-Baxter system can be decomposed as a direct sum of two semigroups. A factorization theorem is proved, generalizing the factorization theorem of Rota-Baxter groups. The notion of Rota-Baxter systems of Lie algebras was introduced, as a generalization of Rota-Baxter Lie algebras. The connection between Rota-Baxter systems of Lie algebras and Lie groups is studied. Finally, as a generalization of the modified Yang-Baxter equation, we define twisted modified Yang-Baxter equations. We give solutions of twisted modified Yang-Baxter equations by Rota-Baxter systems of Lie algebras.
\end{abstract}

{\small

}

%%---------------------------------------------------------------------------
%%------------------------Content-------------------------------------------
%%----------------------------------------------------------------------------

\section{Introduction}\label{Sec:Intro}

\newcommand{\RNum}[1]{\uppercase\expandafter{romannumeral #1\relex}}

The quantum Yang-Baxter equation arose from the study of statistical mechanics and quantum field theory \cite{Baxter,Baxter1,Yang}. In \cite{Drinfeld}, V. G. Drinfeld suggested to study set-theoretical solutions to such equations. The classical Yang-Baxter equation was recognised as the semi-classical limit of the quantum Yang-Baxter equation in the 1980s. Its connection with classical integrable systems and quantum groups was studied in \cite{Chari}.

In 1982, A. A. Belavin and V. G. Drinfeld studied nondegenerate solutions of the classical Yang-Baxter equation for a simple Lie algebra in their remarkable paper \cite{Belavin}. In \cite{Semenov}, M. A. Semenov-Tian-Shansky showed that if there exists a nondegenerate symmetric invariant bilinear form on a Lie algebra $(\mathfrak{g},[\cdot,\cdot])$, then a skew-symmetric solution of the classical Yang-Baxter equation on $\mathfrak{g}$ can be expressed as a linear operator $B: \mathfrak{g} \to \mathfrak{g}$ satisfying
$$[B(u),B(v)]=B([B(u),v])+B([u,B(v)]),\quad \forall u,v\in \mathfrak{g}.$$
The modified Yang-Baxter equation
\begin{equation}\label{MYB}
[R(u),R(v)]=R([R(u),v])+R([u,R(v)])-[u,v],\quad \forall u,v\in \mathfrak{g}
\end{equation}
was also introduced in 1983 in \cite{Semenov}. The infinitesimal factorization theorem for the Lie algebra $\mathfrak{g}$ was obtained from the modified Yang-Baxter equation and the global factorization theorem for the corresponding Lie group was obtained after integration. These theorems have led to important applications in theoretical physics and statistical mechanics \cite{Frenkel,Izosimov,Li,Reshetikhin1,Reshetikhin2}.

In 1982, the concept of Rota-Baxter operators of weight $-1$ appeared in  \cite{Belavin}. More generally, the Rota-Baxter operators of weight $\lambda$ were studied in \cite{Semenov}.
A linear operator $B:\mathfrak{g}\to\mathfrak{g}$ on a Lie algebra $(\mathfrak{g},[\cdot,\cdot])$ is called a Rota-Baxter operator of weight $\lambda$ if
$$[B(u),B(v)]=B([B(u),v])+B([u,B(v)])+\lambda B([u,v]), \quad \forall u,v\in\mathfrak{g}.$$
A Lie algebra with a Rota-Baxter operator of weight $\lambda$ is called a Rota-Baxter Lie algebra of weight $\lambda$. Note that if $R=\operatorname{id}+2B$, then $R:\mathfrak{g}\to\mathfrak{g}$ satisfies the modified Yang-Baxter equation if and only if $B:\mathfrak{g}\to\mathfrak{g}$ is a Rota-Baxter operator of weight $1$.

Recall that a set-theoretical solution to the quantum Yang-Baxter equation is a pair $(X,r)$, where $X$ is a set and $r:X\times X\to X\times X$ is a bijective map such that
$$(r\times \operatorname{id})(\operatorname{id}\times r)(r\times \operatorname{id})=(\operatorname{id}\times r)(r\times \operatorname{id})(\operatorname{id}\times r).$$
In 2006, braces were introduced by W. Rump in \cite{Rump,Rump1} as tools to study involutive set-theoretical solutions to the quantum Yang-Baxter equation. In 2014, F. Ced\'o, E. Jespers and J. Okni\'nski defined left braces in \cite{Cedo}, which are  equivalent to the braces defined by W. Rump. In 2017, L. Guarnieri and L. Vendramin defined skew left braces \cite{LG}, which give non-involutive solutions to the quantum Yang-Baxter equation. Recall that a set $A$ with two binary operations $\cdot$ and $\circ$, with each making $A$ into a group, is called a skew left brace if for all $a,b,c\in A$,
$$a\circ(b\cdot c)=(a\circ b)\cdot a^{-1}\cdot(a\circ c).$$
If $(G,\cdot)$ is abelian, then the  skew left brace $(G,\cdot,\circ)$ becomes a left brace.

In 2021, the Rota-Baxter Lie groups were introduced by L. Guo, H. Lang and Y. Sheng in \cite{LG1} as the integration of Rota-Baxter Lie algebras of weight $1$. Recall that a Rota-Baxter group is a group $G$ with an operator $\mathcal{B}:G\to G$ such that
$$\mathcal{B}(a)\mathcal{B}(b)=\mathcal{B}(a\mathcal{B}(a)b\mathcal{B}(a)^{-1}),\quad \forall a,b\in G.$$
With the help of Rota-Baxter operators, the factorization theorem of groups in \cite{Semenov} can be directly achieved.

In 2022, V. G. Bardakov and V. Gubarev studied the deep connection between skew left braces and Rota-Baxter groups in \cite{VG}. It was shown that there is a skew left brace structure on every Rota-Baxter group and every skew left brace can be injectively embedded into a Rota-Baxter group.

In an attempt to understand the origins and the nature of the law binding two group operations together into skew braces, the notion of skew trusses was proposed by T. Brzezi\'nski in \cite{TB2} in 2019. A skew truss consists of a set $G$ with a group operation $\cdot$ and a semigroup operation $\circ$ connected by a modified distributive law such that for any $a,b,c\in G$
$$a\circ(b\cdot c)=(a\circ b)\cdot \Phi(a)^{-1}\cdot(a\circ c),$$
where $\Phi:G\to G$ is an operator on $G$ called cocycle. Especially, when $\Phi=\operatorname{id}$, the skew truss $(G,\cdot,\circ)$ is a skew left brace.

In light of the fact every Rota-Baxter group has the structure of skew braces, it is natural to address the following questions:
 \begin{description}
 \item[(i)] are there analogues of Rota-Baxter groups which have the structure of skew trusses?
 \item[(ii)] what is the connection between such algebraic systems and Rota-Baxter groups?
 \item[(iii)] whether there is a factorization theorem on such algebraic systems or not?
 \item[(iv)] what is the corresponding structure on the Lie algebras of such algebraic systems?
 \item[(v)] can such algebraic systems give solutions to the Yang-Baxter equation or the modified Yang-Baxter equation?
\end{description}

To answer the first question, we introduce the notion of Rota-Baxter systems of (semi)groups. A Rota-Baxter system is a triple $(G,\mathcal{B}_1,\mathcal{B}_2)$ with $(G,\cdot)$ a group, $\mathcal{B}_1:G\to G$ and $\mathcal{B}_2:G\to G$ two operators on $G$ satisfying \eqref{RB1} and \eqref{RB2}.  It is shown that Rota-Baxter systems are natural generalizations of Rota-Baxter groups. Moreover, $\mathcal{B}_1$ and $\mathcal{B}_2$ induce a semigroup multiplication $\circ$ on $G$, called the descendent operation. We obtain that $(G,\cdot,\circ)$ is a skew truss. The cocycle map $\Phi:G\to G$ is define by $\Phi(a)=\mathcal{B}_1(a)\mathcal{B}_2(a)$ for any $a\in G$. We also show that $G$ has the structure of Rota-Baxter groups if $(G,\mathcal{B}_1,\mathcal{B}_2)$ is a Rota-Baxter system of groups and $\Phi$ is bijective on $G$, which answers the second question.

Similar to the decomposition theorem of Rota-Baxter groups, we obtain a decomposition theorem on Rota-Baxter systems, which answers the third question. In order to answer the fourth question, we introduce Rota-Baxter systems of Lie algebras and prove that if there is a Rota-Baxter system structure on a Lie group, then there is a Rota-Baxter system structure on its Lie algebra under certain conditions. To answer the last question, we define the twisted modified Yang-Baxter equations and obtain that two operators on a Lie algebra satisfy the twisted modified Yang-Baxter equations, if and only if the Lie algebra is a Rota-Baxter system of Lie algebras.

We give the context of the paper. In Section \ref{2.0}, we first introduce the concept of Rota-Baxter systems. As a special example of Rota-Baxter systems, the notion of a twisted Rota-Baxter operator is introduced. In Section \ref{3}, we study the relation between Rota-Baxter systems and skew trusses. The connection between Rota-Baxter systems and Rota-Baxter groups is also considered. In Section \ref{4}, we study the structure of Rota-Baxter systems. In Section \ref{5}, a factorization theorem of a Rota-Baxter system is derived. In Section \ref{6}, the notion of the twisted modified Yang-Baxter equations is proposed. We study the connection between the solutions of twisted modified Yang-Baxter equations and a Rota-Baxter system.

Throughout the paper, for a monoid $G$, we denote the identity element of $G$ by $1_G$.

%%%%%%%%%%%%%%%%%%%%%%%%%%%%%%%%%%%%%%%%%%%%%% %%%%%%%%%%%%%%%%%%%%%%%%%%%%%%%%%%%%%%%%%%%%%%%%%%%%%%%%%%%%%%%%%%%%%%%%%

\section{Rota-Baxter systems}\label{2.0}

In this section, we define Rota-Baxter systems of (semi)groups and give some examples.

We recall the definition of Rota-Baxter groups from \cite{LG1}. Let $G$ be a group. A map $\mathcal{B}:G\to G$ such that
$$\mathcal{B}(a)\mathcal{B}(b)=\mathcal{B}(a\mathcal{B}(a)b\mathcal{B}(a)^{-1}),\quad \forall a,b\in G$$
is called a Rota-Baxter operator of weight $1.$ A map $\mathcal{C}:G\to G$ such that
$$\mathcal{C}(a)\mathcal{C}(b)=\mathcal{C}(\mathcal{C}(a)b\mathcal{C}(a)^{-1}a),\quad \forall a,b\in G$$
is called a Rota-Baxter operator of weight $-1.$ Without further mention, we will call Rota-Baxter operators of weight $1$ simply Rota-Baxter operators. A group with a Rota-Baxter operator is called a Rota-Baxter group.

As a generation of Rota-Baxter groups, we introduce the notion of Rota-Baxter systems.

\begin{Def}\label{RBS}
A Rota-Baxter system of semigroups is a triple $(G,\mathcal{B}_1,\mathcal{B}_2)$ with a semigroup $G$ and two operators $\mathcal{B}_{1}:G\to G$ and $\mathcal{B}_{2}:G\to G$  such that for any $a,b\in G$,
\begin{subequations}
\begin{align}
\mathcal{B}_1(a)\mathcal{B}_1(b)=\mathcal{B}_1(\mathcal{B}_1(a)b\mathcal{B}_2(a)),\label{RB1}\\
\mathcal{B}_{2}(b)\mathcal{B}_{2}(a)=\mathcal{B}_{2}(\mathcal{B}_{1}(a)b\mathcal{B}_{2}(a)).\label{RB2}
\end{align}
\end{subequations}
Moreover, if $G$ is a group, we call $(G,\mathcal{B}_1,\mathcal{B}_2)$ a Rota-Baxter system of groups.
\end{Def}

Note that \eqref{RB1} and \eqref{RB2} imply that both $\mathcal{B}_1(G)$ and $\mathcal{B}_2(G)$ are subsemigroups of $G$.

As explained in the following lemma, Rota-Baxter groups are examples of Rota-Baxter systems of groups.

\begin{lem}
Let $G$ be a Rota-Baxter group with the Rota-Baxter operator $\mathcal{B}$. Define two operators $\mathcal{B}_1:G\to G$ and $\mathcal{B}_2:G\to G$ by
$$\mathcal{B}_{1}(a)=a\mathcal{B}(a), \quad \mathcal{B}_{2}(a)=\mathcal{B}(a)^{-1}, \quad \forall a\in G.$$
Then $(G, \mathcal{B}_{1},\mathcal{B}_{2})$ is a Rota-Baxter system of groups.
\end{lem}

\proof
By the proof of \cite[Proposition 3.1]{LG1}, \eqref{RB1} holds. It remains to show \eqref{RB2}. In fact, for any $a, b\in G$, we have
\begin{align*}
\mathcal{B}_{2}(b)\mathcal{B}_{2}(a)
=&\mathcal{B}(b)^{-1}\mathcal{B}(a)^{-1}
=(\mathcal{B}(a)\mathcal{B}(b))^{-1}\\
=&\mathcal{B}(a\mathcal{B}(a)b\mathcal{B}(a)^{-1})^{-1}
=\mathcal{B}_2(a\mathcal{B}(a)b\mathcal{B}(a)^{-1})\\
=&\mathcal{B}_{2}(\mathcal{B}_1(a)b\mathcal{B}_2(a)),
\end{align*}
as desired.
\qed

The following result shows how to construct Rota-Baxter systems from two subsemigroups.

\begin{prop}\label{RGB4}
Let $G$ be a semigroup and  let $G_1, G_2$ be two subsemigroups of $G$ such that every element $a$ of $G$ can be decomposed uniquely as
$$a=a_1a_2, \qquad a_1\in G_1,\;a_2\in G_2.$$
Define $\mathcal{B}_{1}:G\to G$ and $\mathcal{B}_2:G\to G$ by
$$\mathcal{B}_{1}(a)=a_1, \qquad \mathcal{B}_2(a)=a_2,$$
Then $(G,\mathcal{B}_1,\mathcal{B}_2)$ is a Rota-Baxter system of semigroups.
\end{prop}

\proof
For any $a=a_1a_2, b=b_1b_2\in G$ with $a_1,b_1\in G_1$ and $a_2,b_2\in G_2$, we have $\mathcal{B}_{1}(a)b\mathcal{B}_2(a)=a_1b_1b_2a_2$. As $a_1b_1\in G_1$ and $b_2a_2\in G_2$, we get
$$\mathcal{B}_1(\mathcal{B}_{1}(a)b\mathcal{B}_2(a))=a_1b_1=\mathcal{B}_1(a)\mathcal{B}_1(b)$$
and
$$\mathcal{B}_2(\mathcal{B}_{1}(a)b\mathcal{B}_2(a))=b_2a_2=\mathcal{B}_2(b)\mathcal{B}_2(a),$$
which show that $(G,\mathcal{B}_1,\mathcal{B}_2)$ is a Rota-Baxter system of semigroups.
\qed

As a specific example, let $G$ be a group and $G_{+},G_{-}$ be subgroups of $G$ such that $G=G_{+}G_{-}$ and $G_{+}\cap G_{-}=\{1_G\}.$ It is easy to verify that any $a\in G$ can be written in a unique way as $a=a_+a_-$ with $a_+\in G_+$ and $a_-\in G_-$. Then $(G,\mathcal{B}_1,\mathcal{B}_2)$ is a Rota-Baxter system of groups with $\mathcal{B}_{1},\mathcal{B}_2:G\to G$ defined as
$$\mathcal{B}_1(a)=a_+, \quad\mathcal{B}_2(a)=a_-,\quad \forall a=a_+a_-,\;a_+\in G_+,a_-\in G_-.$$
We mention that  as proved in \cite[Lemma 2.6]{LG1}, the operator $\mathcal{B}:G\to G$ defined by $\mathcal{B}(a)=\mathcal{B}_2(a)^{-1}$ is a Rota-Baxter operator on $G$.

For another example, we consider the monoid on the set of two elements.

\begin{eg}
Let $G$ be the free monoid on the set $\{x,y\}$, whose elements are all the finite strings of zero or more elements from $\{x,y\}$. The multiplication of the monoid is given by string concatenation, and the identity $1$ of the monoid is the unique string of zero elements. Let $G_1$ be the free monoid on the set $\{x\}$, and let $G_2$ be the subset of $G$ which contains the identity $1$  and the strings beginning with $y$. Obviously, $G_1$ and $G_2$ are subsemigroups of $G$, and any element $a\in G$ can be decomposed uniquely as $a=a_1a_2$ with  $a_1\in G_1$ and $a_2\in G_2$. Then we get a Rota-Baxter system of semigroups $(G,\mathcal{B}_1,\mathcal{B}_2)$, where $\mathcal{B}_1$ and $\mathcal{B}_2$ are defined as in Proposition \ref{RGB4}.
\end{eg}

The next example shows how to get a Rota-Baxter system structure from the direct product of two monoids.

\begin{eg}
Let $G$ and $H$ be two monoids. Let $\mathcal{B},\mathcal{L}:G\times H\to G\times H$ be the operators defined by
$$\mathcal{B}_1((a,b))=(a,1_H),\quad \mathcal{B}_2((a,b))=(1_G,1_H),\quad \forall a \in G,\;\forall b\in H.$$
Then $(G\times H,\mathcal{B}_1,\mathcal{B}_2)$ is a Rota-Baxter system of semigroups.
\end{eg}

To give more examples of Rota-Baxter systems, here we propose the notion of twisted Rota-Baxter operators, which is the group version of the twisted Rota-Baxter operators defined in \cite{TB1}.

\begin{Def}
Let $G$ be a group and let $\mathcal{B},\Psi:G\to G$ be two operators with $\Psi$ a homomorphism of groups. The operator $\mathcal{B}$ is called a $\Psi$-twisted Rota-Baxter operator if for any $a, b\in G$,
$$\mathcal{B}(a)\mathcal{B}(b)=\mathcal{B}(\mathcal{B}(a)b\Psi(\mathcal{B}(a)^{-1})).$$
\end{Def}

\begin{prop}\label{TWS}
Let $G$ be a group and $\Psi:G\to G$ be a homomorphism of groups. For a $\Psi$-twisted operator $\mathcal{B}:G\to G$, define $\mathcal{B}_{1}=\mathcal{B}$ and $\mathcal{B}_{2}:G\to G$ such that $\mathcal{B}_2(a)=\Psi(\mathcal{B}(a)^{-1})$ for any $a\in G$. Then $(G, \mathcal{B}_1,\mathcal{B}_2)$ is a Rota-Baxter system of groups.
\end{prop}

\proof
One can readily verify \eqref{RB1} from the definition of the twisted Rota-Baxter operators. It remains to show \eqref{RB2}. For any $a, b\in G$, we have
\begin{align*}
\mathcal{B}_{2}(b)\mathcal{B}_{2}(a)
=&\Psi(\mathcal{B}(b)^{-1})\Psi(\mathcal{B}(a)^{-1})
=\Psi(\mathcal{B}(b)^{-1}\mathcal{B}(a)^{-1})\\
=&\Psi((\mathcal{B}(a)\mathcal{B}(b))^{-1})
=\Psi(\mathcal{B}(\mathcal{B}(a)b\Psi(\mathcal{B}(a)^{-1}))^{-1})\\
=&\mathcal{B}_2(\mathcal{B}_1(a)b\mathcal{B}_2(a)),
\end{align*}
as required.
\qed

Using  Proposition \ref{TWS}, we have the following specific example of Rota-Baxter systems of groups.

\begin{eg}\label{TWS1}
Let $G$ be a group and $\mathcal{B}_1:G\to G$ be an operator. Define $\Psi:G\to G$ by $\Psi(a)=1_G$ for any $a\in G$. Then $\mathcal{B}_1$ is a $\Psi$-twisted Rota-Baxter operator if and only if
$$\mathcal{B}_1(a)\mathcal{B}_1(b)=\mathcal{B}_1(\mathcal{B}_1(a)b),\quad \forall a,b\in G.$$
Assume that $\mathcal{B}_1$ is a $\Psi$-twisted Rota-Baxter operator, then $(G,\mathcal{B}_1,\mathcal{B}_2)$ is a Rota-Baxter system of groups with $\mathcal{B}_2=\Psi$.
\end{eg}

%%%%%%%%%%%%%%%%%%%%%%%%%%%%%%%%%%%%%%%%%%%%%%%%%%%%%%%%%%%%%%%%%%%%%%%%%%%%%%%%%%%%%%%%%%%%%%%%%%%%%%%%%%%%%%
\section{Skew trusses and descendent operations}\label{3}
%\subsection{}
In this section, the relationship between  Rota-Baxter systems and  skew trusses is studied. We introduce the notion of the descendent operations and study the properties of such operations. Moreover, we give a condition to characterize when a Rota-Baxter system has a group structure with respect to the descendent operation.

From here to the end of the paper, a Rota-Baxter system always means a Rota-Baxter system of groups.

\subsection{Skew trusses and Rota-Baxter systems}

In this subsection, we propose the notion of descendent operations and discuss the relationship between Rota-Baxter systems and skew trusses.

Firstly, let's recall the definition of skew trusses given in \cite{TB2}. A skew left truss is a quadruple $(G,\cdot,\circ,\Phi)$, where $G$ is a group with multiplication $\cdot$, $\circ:G\times G\to G$ is an associative operation, and $\Phi:G\to G$ is an operator satisfying
$$a\circ (b\cdot c)=(a\circ b)\cdot\Phi(a)^{-1}\cdot(a\circ c),\quad \forall a, b, c\in G.$$
The operator $\Phi$ is called a cocycle. Note that the cocycle of a skew left truss is denoted by $\sigma$ in \cite{TB2}. Without further mention, we will call skew left trusses simply by skew trusses in this paper.

Now we introduce the definition of descendent operations.

\begin{Def}\label{CIR}
 Let $G$ be a group with two operators $\mathcal{B}_1:G\to G$ and $\mathcal{B}_2:G\to G$. The operation $\circ:G\times G\to G$ defined by
\begin{equation}\label{CP}
a\circ b=\mathcal{B}_{1}(a)b\mathcal{B}_{2}(a),\quad  \forall a, b\in G
\end{equation}
is called the descendent operation of $(G,\mathcal{B}_1,\mathcal{B}_2)$. The descendent operation $\circ$ is called associative if it satisfies the associative law, that is $a\circ(b\circ c)=(a\circ b)\circ c$ for any $a,b,c\in G.$ If $\circ$ is associative, then the operator $\Phi:G\to G$ given by
 $$\Phi(a)=\mathcal{B}_1(a)\mathcal{B}_2(a)$$
is called the cocycle of the triple $(G,\mathcal{B}_1,\mathcal{B}_2)$.
\end{Def}

The following proposition gives equivalent characterisations of a Rota-Baxter system.

\begin{prop}\label{HOM}
Let $G$ be a group with operators $\mathcal{B}_1:G\to G$ and $\mathcal{B}_2:G\to G.$ Let $\circ$ be the descendent operation of $(G,\mathcal{B}_1,\mathcal{B}_2)$. Then the following statements are equivalent:
 \begin{description}
\item[(i)] $(G,\mathcal{B}_1,\mathcal{B}_2)$ is a Rota-Baxter system;
\item[(ii)] $\circ$ is associative and $\mathcal{B}_1$ is a semigroup homomorphism from $(G,\circ)$ to $(G,\cdot)$;
\item[(iii)] $\circ$ is associative and $\mathcal{B}_2$ is a semigroup anti-homomorphism from $(G,\circ)$ to $(G,\cdot)$.
\end{description}
\end{prop}

\proof
For any $a, b, c\in G$,
$$(a\circ b)\circ c=a\circ(b\circ c) \Longleftrightarrow \mathcal{B}_1(a\circ b)c\mathcal{B}_2(a\circ b)
=\mathcal{B}_1(a)\mathcal{B}_1(b)c\mathcal{B}_2(b)\mathcal{B}_2(a).$$

Suppose that (i) holds. Then using the facts $\mathcal{B}_1(a\circ b)=\mathcal{B}_1(a)\mathcal{B}_1(b)$ and $\mathcal{B}_2(a\circ b)=\mathcal{B}_2(b)\mathcal{B}_2(a)$, we find $\circ$ is associative. Hence we get (ii) and (iii).

Suppose that (ii) holds. Since $\circ$ is associative and $\mathcal{B}_1(a\circ b)=\mathcal{B}_1(a)\mathcal{B}_1(b)$, we have
$$\mathcal{B}_1(a\circ b)c\mathcal{B}_2(a\circ b)
=\mathcal{B}_1(a\circ b)c\mathcal{B}_2(b)\mathcal{B}_2(a).$$
Then $\mathcal{B}_2(a\circ b)=\mathcal{B}_2(b)\mathcal{B}_2(a)$ and (iii) holds. Similarly, we can get (ii) from (iii).

Now it is obvious that (ii) (resp. (iii)) implies (i). \qed

Next we show that every Rota-Baxter system is a skew truss.

\begin{prop}\label{TRU}
Let $G$ be a group with two operators $\mathcal{B}_1:G\to G$ and $\mathcal{B}_2:G\to G$. If the descendent operation $\circ$ of $(G,\mathcal{B}_1,\mathcal{B}_2)$ is associative, then $(G,\cdot,\circ,\Phi)$ is a skew truss, where $\Phi$ is the cocycle of $(G,\mathcal{B}_1,\mathcal{B}_2)$. Especially, every Rota-Baxter system is a skew truss.
\end{prop}

\proof
For any $a,b,c\in G$, we have
\begin{align*}
a\circ (bc)=&\mathcal{B}_1(a)(bc)\mathcal{B}_2(a)=(\mathcal{B}_1(a)b\mathcal{B}_2(a))\Phi (a)^{-1}(\mathcal{B}_1(a)c\mathcal{B}_2(a))\\
=&(a\circ b)\Phi (a)^{-1}(a\circ c),
\end{align*}
as desired.
\qed

The following lemma will be used later.

\begin{lem}\label{TB}
Let $G$ be a group with two operators $\mathcal{B}_1:G\to G$ and $\mathcal{B}_2:G\to G$. If the descendent operation $\circ$ of $(G,\mathcal{B}_1,\mathcal{B}_2)$ is associative, then the cocycle $\Phi$ of $(G,\mathcal{B}_1,\mathcal{B}_2)$ satisfies $\Phi(a\circ b)=a\circ \Phi(b)$ for any $a,b\in G$.
\end{lem}

\proof
The result follows from Proposition \ref{TRU} and \cite[Lemma 2.3]{TB2}.
\qed

 By \cite[Theorem 2.9]{TB2}, $(G,\circ)$ acts from the left on $(G,\cdot)$  if $(G,\cdot,\circ)$ is a skew truss. For a triple $(G,\mathcal{B}_1,\mathcal{B}_2)$ where $G$ is a group with two operators $\mathcal{B}_1:G\to$ and $\mathcal{B}_2:G\to G$, we show that such action exists if and only if only if $(G,\mathcal{B}_1,\mathcal{B}_2)$ has the structure of skew trusses.

\begin{thm}\label{TH1}
Let $G$ be a group with two operators $\mathcal{B}_1:G\to G$ and $\mathcal{B}_2:G\to G$. Let $\circ$ be the descendent operation and let the maps $\lambda:G\to \operatorname{Aut}(G);a\mapsto \lambda_a$ and $\mu:G\to \operatorname{Aut}(G);a\mapsto \mu_a$ be defined as
$$\lambda_{a}(b)=\mathcal{B}_2(a)^{-1}b\mathcal{B}_2(a),\quad \mu_{a}(b)=\mathcal{B}_1(a)b\mathcal{B}_1(a)^{-1},\quad \forall b\in G,$$
respectively. Then the following statements are equivalent:
 \begin{description}
\item[(i)] $\circ$ is associative and then $(G,\cdot,\circ,\Phi)$ is a skew truss;
\item[(ii)] $\lambda_a\lambda_b=\lambda_{a\circ b}$ and $(a\circ b)\circ 1_G=a\circ(b\circ 1_G)$ for any $a,b\in G$;
\item[(iii)] $\mu_a\mu_b=\mu_{a\circ b}$ and $(a\circ b)\circ 1_G=a\circ(b\circ 1_G)$ for any $a,b\in G$.
\end{description}
\end{thm}

\proof
Assume (i) holds. For any $a,b,c\in G$,
\begin{align*}
\lambda_a\lambda_b(c)
&=\mathcal{B}_2(a)^{-1}\mathcal{B}_2(b)^{-1}c\mathcal{B}_2(b)\mathcal{B}_2(a)\\
&=(\mathcal{B}_1(a)\mathcal{B}_1(b)\mathcal{B}_2(b)\mathcal{B}_2(a))^{-1}(\mathcal{B}_1(a)\mathcal{B}_1(b)c\mathcal{B}_2(b)\mathcal{B}_2(a))\\
&=(a\circ (b\circ 1_G))^{-1}(a\circ (b\circ c)).
\end{align*}
As $\circ$ is associative, we get
\begin{align*}
\lambda_a\lambda_b(c)
=&((a\circ b)\circ 1_G)^{-1}((a\circ b)\circ c)\\
=&(\mathcal{B}_1(a\circ b)\mathcal{B}_2(a\circ b))^{-1}(\mathcal{B}_1(a\circ b)c\mathcal{B}_2(a\circ b))\\
=&\mathcal{B}_2(a\circ b)^{-1}c\mathcal{B}_2(a\circ b)=\lambda_{a\circ b}(c).
\end{align*}
Similarly, we have
\begin{align*}
\mu_a\mu_b(c)
&=\mathcal{B}_1(a)\mathcal{B}_1(b)c\mathcal{B}_1(b)^{-1}\mathcal{B}_1(a)^{-1}\\
&=(\mathcal{B}_1(a)\mathcal{B}_1(b)c\mathcal{B}_2(b)\mathcal{B}_2(a))(\mathcal{B}_1(a)\mathcal{B}_1(b)\mathcal{B}_2(b)\mathcal{B}_2(a))^{-1}\\
&=(a\circ(b\circ c))(a\circ (b\circ 1_G))^{-1}=((a\circ b)\circ c)((a\circ b)\circ 1_G)^{-1}\\
&=(\mathcal{B}_1(a\circ b)c\mathcal{B}_2(a\circ b))(\mathcal{B}_1(a\circ b)\mathcal{B}_2(a\circ b))^{-1}\\
&=\mathcal{B}_1(a\circ b)c\mathcal{B}_1(a\circ b)^{-1}=\mu_{a\circ b}(c).
\end{align*}
Hence $\lambda_a\lambda_b=\lambda_{a\circ b}$ and $\mu_a\mu_b=\mu_{a\circ b}$. And then we get (ii) and (iii).

Conversely, assume (ii) holds. For any $a,b,c\in G$,
\begin{align*}
a\circ(b\circ c)
&=a\circ(\Phi(b)\mathcal{B}_2(b)^{-1}c\mathcal{B}_2(b))
=a\circ(\Phi(b)\lambda_b(c))\\
&=\mathcal{B}_1(a)\Phi(b)\lambda_b(c)\mathcal{B}_2(a)
=\Phi(a)\lambda_a(\Phi(b)\lambda_b(c))\\
&=\Phi(a)\lambda_a(\Phi(b))\lambda_a\lambda_b(c)
=\mathcal{B}_1(a)\mathcal{B}_1(b)\mathcal{B}_2(b)\mathcal{B}_2(a)\lambda_a\lambda_b(c).
\end{align*}
As $\lambda_a\lambda_b=\lambda_{a\circ b}$ and $(a\circ b)\circ 1_G=a\circ(b\circ 1_G)$, we get
\begin{align*}
a\circ(b\circ c)
&=(a\circ(b\circ 1_G))\lambda_{a\circ b}(c)
=((a\circ b)\circ 1_G)\lambda_{a\circ b}(c)\\
&=(\mathcal{B}_1(a\circ b)\mathcal{B}_2(a\circ b))\lambda_{a\circ b}(c)
=\mathcal{B}_1(a\circ b)c\mathcal{B}_2(a\circ b)\\
&=(a\circ b)\circ c,
\end{align*}
which implies (i).

Similarly, if (iii) holds, then for any $a,b,c\in G$, we have
\begin{align*}
a\circ(b\circ c)
&=a\circ(\mathcal{B}_1(b)c\mathcal{B}_1(b)^{-1}\Phi(b))
=a\circ(\mu_b(c)\Phi(b))\\
&=\mathcal{B}_1(a)\mu_b(c)\Phi(b)\mathcal{B}_2(a)
=\mu_a(\mu_b(c)\Phi(b))\Phi(a)\\
&=\mu_a\mu_b(c)\mu_a(\Phi(b))\Phi(a)
=\mu_a\mu_b(c)\mathcal{B}_1(a)\mathcal{B}_1(b)\mathcal{B}_2(b)\mathcal{B}_2(a)\\
&=\mu_a\mu_b(c)(a\circ(b\circ 1_G))
=\mu_a\circ_b(c)((a\circ b)\circ 1_G)\\
&=\mu_{a\circ b}(c)(\mathcal{B}_1(a\circ b)\mathcal{B}_2(a\circ b))
=\mathcal{B}_1(a\circ b)c\mathcal{B}_2(a\circ b)\\
&=(a\circ b)\circ c,
\end{align*}
which deduces (i).\qed

As a consequence of the above theorem, we get

\begin{cor}\label{CEN}
Let $G$ be a group with operators $\mathcal{B}_1:G\to G$ and $\mathcal{B}_2:G\to G.$ If the descendent operation $\circ$ of $(G,\mathcal{B}_1,\mathcal{B}_2)$ is associative and the center of $G$ is trivial, then $(G,\mathcal{B}_1,\mathcal{B}_2)$ is a Rota-Baxter system.
\end{cor}

\proof
If $\circ$ is associative, then from Theorem \ref{TH1}, for any $a,b,c\in G$, we have
$$\mathcal{B}_2(a)^{-1}\mathcal{B}_2(b)^{-1}c\mathcal{B}_2(b)\mathcal{B}_2(a)=\mathcal{B}_2(a\circ b)^{-1}c\mathcal{B}_2(a\circ b).$$
Hence
$$c\mathcal{B}_2(b)\mathcal{B}_2(a)\mathcal{B}_2(a\circ b)^{-1}=\mathcal{B}_2(b)\mathcal{B}_2(a)\mathcal{B}_2(a\circ b)^{-1}c,$$
which means $\mathcal{B}_2(b)\mathcal{B}_2(a)\mathcal{B}_2(a\circ b)^{-1}\in Z(G)$, the center of $G$. As $Z(G)$ is trivial, we get
$$\mathcal{B}_2(b)\mathcal{B}_2(a)\mathcal{B}_2(a\circ b)^{-1}=1_G,$$
which is equivalent to
$$\mathcal{B}_2(b)\mathcal{B}_2(a)=\mathcal{B}_2(a\circ b).$$
Similarly, we can show
$$\mathcal{B}_1(a)\mathcal{B}_1(b)=\mathcal{B}_1(a\circ b).$$
Hence $(G,\mathcal{B}_1,\mathcal{B}_2)$ is a Rota-Baxter system.
\qed

\subsection{Left identities and right inverses in $(G,\circ)$}

Let $G$ be a group with operators $\mathcal{B}_1:G\to G$ and $\mathcal{B}_2:G\to G.$ For any $a\in G$, denote  $\mathcal{B}_1(a)^{-1}a\mathcal{B}_2(a)^{-1}$ by $e_a$.
We show that if the descendent operation $\circ$ of $(G,\mathcal{B}_1,\mathcal{B}_2)$ is associative, then $e_a$ ($a\in G$) are left identities of $(G,\circ)$.

\begin{thm}\label{SEM}
Let $G$ be a group with operators $\mathcal{B}_1:G\to G$ and $\mathcal{B}_2:G\to G.$ If the descendent operation $\circ$ of $(G,\mathcal{B}_1,\mathcal{B}_2)$ is associative, then for any $a\in G$, $e_a$ is a left identity of the semigroup $(G,\circ)$. Moreover, if $(G,\mathcal{B}_1,\mathcal{B}_2)$ is a Rota-Baxter system, then
$$\mathcal{B}_1(e_a)=\mathcal{B}_2(e_a)=1_G.$$
\end{thm}

\proof
For any $a\in G$, we have $a\circ e_a=a$. Hence from the proof of Corollary \ref{CEN}, we find
$$\mathcal{B}_2(e_a)=\mathcal{B}_2(e_a)\mathcal{B}_2(a)\mathcal{B}_2(a)^{-1}\in Z(G).$$
Similarly, we have
$$\mathcal{B}_1(e_a)=\mathcal{B}_1(a)^{-1}\mathcal{B}_1(a)\mathcal{B}_1(e_a)\in Z(G).$$
As $\circ$ is associative, we have $a\circ (e_a\circ 1_G)=(a\circ e_a)\circ 1_G$, which is
$$\mathcal{B}_1(a)\mathcal{B}_1(e_a)\mathcal{B}_2(e_a)\mathcal{B}_2(a)=\mathcal{B}_1(a)\mathcal{B}_2(a).$$
Therefore we obtain
\begin{gather}\label{SIGMA}
\mathcal{B}_1(e_a)\mathcal{B}_2(e_a)=1_G.
\end{gather}
Hence for any $b\in G$, we have
$$e_a\circ b= \mathcal{B}_1(e_a)b\mathcal{B}_2(e_a)
=b\mathcal{B}_1(e_a)\mathcal{B}_2(e_a)=b,$$
which means that $e_a$ is a left identity of the semigroup $(G,\circ)$.

When $(G,\mathcal{B}_1,\mathcal{B}_2)$ is a Rota-Baxter system, we have
$$\mathcal{B}_1(a)\mathcal{B}_1(e_a)=\mathcal{B}_1(a), \quad \forall a\in G.$$
Hence $\mathcal{B}_1(e_a)=1_G$ and we have $\mathcal{B}_2(e_a)=1_G$ similarly.
\qed

The following proposition shows that if the descendent operation $\circ$ of $(G,\mathcal{B}_1,\mathcal{B}_2)$ is associative, then $(G,\circ)$ is right divisible.

\begin{prop}\label{INV}
If the descendent operation $\circ$ of $(G,\mathcal{B}_1,\mathcal{B}_2)$  is associative, then $(G,\circ)$ is right divisible: for any $a,b\in G$, there exists a unique $c\in G$ such that
$$a\circ c=b.$$
\end{prop}

\proof
Taking $c=\mathcal{B}_1(a)^{-1}b\mathcal{B}_2(a)^{-1}$, we see $a\circ c=b$. The uniqueness of $c$ follows from the fact that $e_a$ is a left identity of $(G,\circ)$.
\qed

\begin{rem}\label{AUT}
Let $(G,\cdot,\circ,\Phi)$ be a skew truss. For any $a\in G$, define the operator $\lambda_a:G\to G$ by
$$\lambda_a(b)=\Phi(a)^{-1}(a\circ b),\qquad \forall b\in G.$$
By \cite[Theorem 2.9]{TB2}, $a\mapsto \lambda_a$ gives an action of the semigroup $(G,\circ)$ on $(G,\cdot)$. We can show that the following statements are equivalent:
\begin{description}
\item[(i)] For any $a, b\in G$, there is a unique $c\in G$ such that $a\circ c=b$;
\item[(ii)] For any $a\in G$, the operator $\lambda_a$ is an automorphism of the group $(G,\cdot)$.
\end{description}

In fact, if $(G,\circ)$ is right divisible, then for any $b\in G$, there exists a unique $c\in G$ such that $a\circ c=\Phi(a)b$. Then we have $\lambda_a(c)=b$ and $\lambda_a$ is surjective. And if $\lambda_a(b_1)=\lambda_a(b_2)$, then $a\circ b_1=a\circ b_2$. Hence $b_1=b_2$ and $\lambda_a$ is injective. Conversely, if $\lambda_a$ is bijective,  then there is a unique $c\in G$ such that $\lambda_a(c)=\Phi(a)^{-1}b$. Hence $a\circ c=b$ and $c$ is unique.
\end{rem}

As in the case of Rota-Baxter groups, if $(G,\mathcal{B}_1,\mathcal{B}_2)$ is a Rota-Baxter system, we can show $\mathcal{B}_1(G)$ and $\mathcal{B}_2(G)$ are subgroups of $(G,\cdot)$. For this purpose, we prepare a lemma.

\begin{lem}\label{DAG}
Let $(G,\mathcal{B}_1,\mathcal{B}_2)$ be a Rota-Baxter system. For any $a\in G$, set $a^{\dag}=\mathcal{B}_1(a)^{-1}e_a\mathcal{B}_2(a)^{-1}$. Then we have
$$\mathcal{B}_1(a^{\dag})=\mathcal{B}_1(a)^{-1},\quad \mathcal{B}_2(a^{\dag})=\mathcal{B}_2(a)^{-1},\quad \forall a\in G.$$
\end{lem}

\proof
By Theorem \ref{SEM}, we have
$$\mathcal{B}_1(a)\mathcal{B}_1(a^{\dag})=\mathcal{B}_1(a\circ a^{\dag})=\mathcal{B}_1(e_a)=1_G,$$
which implies that $\mathcal{B}_1(a)^{-1}=\mathcal{B}_1(a^{\dag})$. We obtain $\mathcal{B}_2(a)^{-1}=\mathcal{B}_2(a^{\dag})$ similarly.
\qed

\begin{cor}\label{INV1}
Let $(G,\mathcal{B}_1,\mathcal{B}_2)$ be a Rota-Baxter system. Then $\mathcal{B}_1(G)$ and $\mathcal{B}_2(G)$ are subgroups of $G$.
\end{cor}

\proof
By \eqref{RB1}, $\mathcal{B}_1(G)$ is closed under the multiplication. And Lemma \ref{DAG} implies that $\mathcal{B}_1(G)$ is closed under the inverse map. Therefore $\mathcal{B}_1(G)$ is a subgroup of $G$. We obtain that $\mathcal{B}_2(G)$ is a subgroup of $G$ similarly.
\qed

\subsection{Rota-Baxter systems and Rota-Baxter groups}

In this subsection, we discuss the relationship between Rota-Baxter systems and Rota-Baxter groups.

By  \cite[Propositon 3.1]{VG}, a Rota-Baxter group has the structure of skew brace. Inspired by \cite[Lemma 4.1]{TB2} which states that a skew truss has the structure of skew brace if its cocycle is bijective, we have the following proposition.

\begin{prop}\label{RBG5}
Let $(G,\mathcal{B}_1,\mathcal{B}_2)$ be a Rota-Baxter system. If the cocycle $\Phi$ is bijective, then the operator $\widetilde{\mathcal{B}}:G\to G$ defined by
$$\widetilde{\mathcal{B}}(a)=\mathcal{B}_1(\Phi^{-1}(a)),\quad \forall a\in G$$
is a Rota-Baxter operator of weight $-1$, and the operator $\mathcal{B}:G\to G$ defined by
$$\mathcal{B}(a)=\mathcal{B}_2(\Phi^{-1}(a))^{-1},\quad \forall a\in G$$
is a Rota-Baxter operator of weight $1$.
\end{prop}

\proof
For any $a, b\in G$, we have
\begin{align*}
\widetilde{\mathcal{B}}(a)\widetilde{\mathcal{B}}(b)&=\mathcal{B}_1(\Phi^{-1}(a))\mathcal{B}_1(\Phi^{-1}(b))=\mathcal{B}_1(\Phi^{-1}(a)\circ \Phi^{-1}(b))\\
&=\mathcal{B}_1( \Phi^{-1}(\Phi(\Phi^{-1}(a)\circ \Phi^{-1}(b))).
\end{align*}
Using Lemma \ref{TB}, we get
$$\widetilde{\mathcal{B}}(a)\widetilde{\mathcal{B}}(b)=\widetilde{\mathcal{B}}(\Phi^{-1}(a)\circ b)=\widetilde{\mathcal{B}}(\widetilde{\mathcal{B}}(a)b\mathcal{B}(a)^{-1}).$$
As
$$\widetilde{\mathcal{B}}(a)\mathcal{B}(a)^{-1}=\mathcal{B}_1(\Phi^{-1}(a))\mathcal{B}_2(\Phi^{-1}(a))=\Phi(\Phi^{-1}(a))=a,$$
we have
$$\widetilde{\mathcal{B}}(a)\widetilde{\mathcal{B}}(b)=\widetilde{\mathcal{B}}(\widetilde{\mathcal{B}}(a)b\widetilde{\mathcal{B}}(a)^{-1}a),$$
which shows that $\widetilde{\mathcal{B}}$ is a Rota-Baxter operator of weight $-1$.

Similarly, for any $a,b\in G$, we have
\begin{align*}
\mathcal{B}(a)\mathcal{B}(b)&=\mathcal{B}_2(\Phi^{-1}(a))^{-1}\mathcal{B}_2(\Phi^{-1}(b))^{-1}=\mathcal{B}_2(\Phi^{-1}(a)\circ \Phi^{-1}(b))^{-1}\\
&=\mathcal{B}( \widetilde{\mathcal{B}}(a)b\mathcal{B}(a)^{-1})=\mathcal{B}(a\mathcal{B}(a)b\mathcal{B}(a)^{-1}),
\end{align*}
which shows that $\mathcal{B}$ is a Rota-Baxter operator of weight $1$.
\qed

We show in the next theorem that a Rota-Baxter system has a group structure with respect to the descendent operation if and only if the cocycle is bijective.

\begin{thm}\label{BIJ}
Suppose that $(G,\mathcal{B}_1,\mathcal{B}_2)$ is a Rota-Baxter system with the descendent operation $\circ$ and the cocycle $\Phi$. Then the following statements are equivalent:
\begin{description}
\item[(i)] $(G,\circ)$ is a group;
\item[(ii)] $\Phi$ is bijective.
\end{description}
\end{thm}

\proof
If $(G,\circ)$ is a group, we get $\Phi$ is bijective from \cite[Lemma 4.1]{TB2}.

Conversely, assume that $\Phi$ is a bijection. For any $a\in G$, by \eqref{SIGMA}, we have $\Phi(e_a)=1_G$, where $e_a=\mathcal{B}_1(a)^{-1}a\mathcal{B}_2(a)^{-1}$. Hence $e_a=\Phi^{-1}(1_G)$ is independent of $a$. We denote $e_a$ by $e$.  By Theorem \ref{SEM}, $e$ is a left identity of  $(G,\circ)$. By the proof of \cite[Lemma 4.1]{TB2}, $e$ is a right identity. Hence $e$ is the identity of $(G,\circ)$.

Set $a^{\dag}=\mathcal{B}_1(a)^{-1}e\mathcal{B}_2(a)^{-1}$. Then $a\circ a^{\dag}=e$. Moreover, Lemma \ref{DAG} claims that
$$\mathcal{B}_1(a^{\dag})=\mathcal{B}_1(a)^{-1} ,\quad \mathcal{B}_2(a^{\dag})=\mathcal{B}_2(a)^{-1}.$$
Hence
$$a^{\dag}\circ a=\mathcal{B}_1(a^{\dag})a\mathcal{B}_2(a^{\dag})=\mathcal{B}_1(a)^{-1}a\mathcal{B}_2(a)^{-1}=e.$$
It follows that $a^{\dag}$ is the inverse of $a$ with respect to $\circ$. Hence $(G,\circ)$ is a group.
\qed

\begin{cor}\label{GRO}
Let $(G,\mathcal{B}_1,\mathcal{B}_2)$ be a Rota-Baxter system with the descendent operation $\circ$ and the cocycle $\Phi$. If $\Phi$ is bijective, then $G$ is a skew brace with respect to $\cdot$ and $\bullet$, where the operation $\bullet:G\times G\to G$ is defined by
$$a\bullet b=\mathcal{B}_1(\Phi^{-1}(a))b\mathcal{B}_2(\Phi^{-1}(a)),\quad \forall a,b\in G.$$
\end{cor}

\proof
The result follows from Theorem \ref{BIJ} and \cite[Lemma 4.3]{TB2}.
\qed

\begin{rem}
In fact, Corollary \ref{GRO} can also be derived from \cite[Propositon 3.1]{VG}, Proposition \ref{RBG5} and Theorem \ref{BIJ}. %In other words, in a Rota-Baxter system $(G,\mathcal{B}_1,\mathcal{B}_2)$ with a bijective cocycle, the skew brace on $G$ introduced in \cite[Lemma 4.3]{TB2} is in coincide with the skew brace on $G$ given in \cite[Propositon 3.1]{VG}.
\end{rem}

\section{Structure of Rota-Baxter systems}\label{4}

In this section, we investigate the structure of Rota-Baxter systems. Firstly, for a Rota-Baxter system, we show that there are groups attached to the left identities obtained in Theorem \ref{SEM}. Then it's shown that every Rota-Baxter system is a union of such groups. Finally we show that every Rota-Baxter system is a direct sum of a subgroup and a subsemigroup.

Let $(G,\mathcal{B}_1,\mathcal{B}_2)$ be a Rota-Baxter system with the descendent operation $\circ$. For any $t\in G$, set
$$G_t=\{a\in G\mid a\circ e_t=a\}$$
and
$$G\circ e_t=\{a\circ e_t\mid a\in G\},$$
where $e_t=\mathcal{B}_1(t)^{-1}t\mathcal{B}_2(t)^{-1}$. It is obvious that $t\in G_t$. We prove that $G_t$ is a group with respect to $\circ$.

\begin{lem}\label{GT}
Let $(G,\mathcal{B}_1,\mathcal{B}_2)$ be a Rota-Baxter system with the descendent operation $\circ$. Then for any $t\in G$, $G_t=G\circ e_t$ and $G_t$ is a group with respect to $\circ$.
\end{lem}

\proof
For any $a\in G_t$, we have $a=a\circ e_t\in G\circ e_t$. Conversely, for any $b=b_1\circ e_t\in G\circ e_t$ with $b_1\in G$, we have
$$b\circ e_t=(b_1\circ e_t)\circ e_t=b_1\circ (e_t\circ e_t)=b_1\circ e_t=b.$$
Hence $b\in G_t$ and then $G_t=G\circ e_t.$

Now we show that $G_t$ is a group with respect to $\circ$. As $e_t=e_t\circ e_t$, we see $e_t\in G_t$. And for any $a\in G_t$, it is obvious that
$$a\circ e_t=a=e_t\circ a.$$
For any $a, b\in G_t,$ we have
$$(a\circ b)\circ e_t=a\circ(b\circ e_t)=a\circ b,$$
which implies that $a\circ b\in G_t$.  For any $a\in G_t$, set $a^{\dag}=\mathcal{B}_1(a)^{-1}e_t\mathcal{B}_2(a)^{-1}$ as before.  Then we have $a\circ a^{\dag}=e_t$. And by Theorem \ref{SEM}, we get $$\mathcal{B}_1(a)\mathcal{B}_1(a^{\dag})=1_G,\quad \mathcal{B}_2(a^{\dag})\mathcal{B}_2(a)=1_G.$$
Therefore
$$a^{\dag}\circ a=\mathcal{B}_1(a)^{-1}a\mathcal{B}_2(a)^{-1}=e_a.$$
Now since $a\circ e_t=a\circ e_a=a$, we have $e_t=e_a$ by Proposition \ref{INV}. Hence we get $a^\dag\circ a=e_t$. Moreover, we have
$$a^{\dag}\circ e_t=a^{\dag}\circ (a\circ a^{\dag})=(a^{\dag}\circ a)\circ a^{\dag}=e_t\circ a^{\dag}=a^{\dag},$$
 which implies that $a^{\dag}\in G_t$. Hence $G_t$ is a group with respect to $\circ$.
\qed

Lemma \ref{GT} claims that $(G_t,\circ)$ is a group with the identity $e_t$.  These groups are in fact isomorphic to each other as proved in the following lemma.

\begin{lem}\label{ISO}
Let $(G,\mathcal{B}_1,\mathcal{B}_2)$ be a Rota-Baxter system. Then for any $t_1,t_2\in G,$ the map
$$\Psi_{t_1,t_2}:(G_{t_1},\circ)\to (G_{t_2},\circ),\quad \Psi_{t_1,t_2}(a)=a\circ e_{t_2},\quad \forall a\in G_{t_1},$$
is a group isomorphism.
\end{lem}

\proof
It is obvious that $\Psi_{t_1,t_2}$ is well-defined. For any $a,b\in G_{t_1}$, we have
\begin{align*}
&\quad \Psi_{t_1,t_2}(a)\circ\Psi_{t_1,t_2}(b)
=(a\circ e_{t_2})\circ (b\circ e_{t_2})\\
&=a\circ(e_{t_2}\circ b)\circ e_{t_2}
=(a\circ b)\circ e_{t_2}=\Psi_{t_1,t_2}(a\circ b).
\end{align*}
Hence $\Psi_{t_1,t_2}$ is a group homomorphism.

For any $c\in G_{t_2}$, we have $c\circ e_{t_1}\in G_{t_1}$ and
$$\Psi_{t_1,t_2}(c\circ e_{t_1})=(c\circ e_{t_1})\circ e_{t_2}=c\circ(e_{t_1}\circ e_{t_2})=c\circ e_{t_2}=c.$$
Hence $\Psi_{t_1,t_2}$ is surjective. It remains to show that $\Psi_{t_1,t_2}$ is injective. Let $a\in G_{t_1}$ such that $\Psi_{t_1,t_2}(a)=e_{t_2}$.  Then we have
$$a=a\circ e_{t_1}=a\circ(e_{t_2}\circ e_{t_1})=(a\circ e_{t_2})\circ e_{t_1}=e_{t_2}\circ e_{t_1}=e_{t_1},$$
which implies that $\Psi_{t_1,t_2}$ is injective.
\qed

The next lemma gives another relation between the groups $G_t$, $t\in G$.

\begin{lem}
Let $(G,\mathcal{B}_1,\mathcal{B}_2)$ be a Rota-Baxter system. Then we have $G_{t_1}\circ G_{t_2}=G_{t_2}$, where
$$G_{t_1}\circ G_{t_2}=\{a\circ b\mid a\in G_{t_1}, b\in G_{t_2}\}.$$
\end{lem}

\proof
For any $a\in G_{t_2}$, we have $a=e_{t_1}\circ a\in G_{t_1}\circ G_{t_2}$. Hence $G_{t_2}\subseteq G_{t_1}\circ G_{t_2}$. Conversely, for any $a\in G_{t_1}$ and $b\in G_{t_2}$, we have
$$a\circ b=(a\circ e_{t_1})\circ (b\circ e_{t_2})=a\circ(e_{t_1}\circ b)\circ e_{t_2}=(a\circ b)\circ e_{t_2}.$$
Therefore $a\circ b\in G_{t_2}$, and then  $G_{t_1}\circ G_{t_2}\subseteq G_{t_2}$.
\qed

Denote the restriction of $\Phi$ on $G_t$ by $\Phi_t,$ that is,
$$\Phi_t:G_t\to G,\quad \Phi_t(a)=\Phi(a),\quad \forall a\in G_t.$$
We compute the image of $\Phi_t$ as shown in the following proposition.

\begin{prop}\label{IMG}
Let $(G,\mathcal{B}_1,\mathcal{B}_2)$ be a Rota-Baxter system. Then for any $t\in G$, $\operatorname{Im}(\Phi_t)=\operatorname{Im}(\Phi)=G_{1_G}$. Moreover, the map $\Phi_{1_G}:G_{1_G}\to G_{1_G}$ is bijective.
\end{prop}

\proof
For any $a\in G$, we have
$$\Phi(a)\circ e_{1_G}=(a\circ 1_G)\circ e_{1_G}=a\circ (1_G\circ e_{1_G})=a\circ 1_G=\Phi(a).$$
Hence $\Phi(a)\in G_{1_G}$, and then $\operatorname{Im}(\Phi)\subseteq G_{1_G}$.

For any $t\in G$, it is obvious that $\operatorname{Im}(\Phi_t)\subseteq \operatorname{Im}(\Phi)$. We want to show $G_{1_G}\subseteq \operatorname{Im}(\Phi_t)$. In fact, for any $a\in G_{1_G}$, set $c=a\circ 1_G^{\dag}$, where $1_G^{\dag}$ is the unique element in $G_{1_G}$ such that
$$1^{\dag}_{G}\circ 1_{G}=1_{G}\circ 1^{\dag}_{G}=e_{1_{G}}.$$
We have $c\circ e_t\in G_t$ and
$$\Phi_t(c\circ e_t)=\mathcal{B}_1(c)\mathcal{B}_1(e_t)\mathcal{B}_2(e_t)\mathcal{B}_2(c)=\mathcal{B}_1(c)\mathcal{B}_2(c)=\Phi(c),$$
where we have used Theorem \ref{SEM}. Now
$$\Phi(c)=\mathcal{B}_1(a\circ 1_G^{\dag})\mathcal{B}_2(a\circ 1_G^{\dag})=(a\circ 1_G^{\dag})\circ 1_G=a\circ(1_G^{\dag}\circ 1_{G})=a\circ e_{1_{G}}=a,$$
which implies that $a\in\operatorname{Im}(\Phi_t)$. Hence we get $G_{1_G}\subseteq \operatorname{Im}(\Phi_t)$. Therefore $\operatorname{Im}(\Phi_t)=\operatorname{Im}(\Phi)=G_{1_G}$.

Since $\operatorname{Im}(\Phi_{1_G})=G_{1_G}$, we have a surjective map $\Phi_{1_G}:G_{1_G}\to G_{1_G}$. It remains to show that $\Phi_{1_G}$ is injective. For any $a_1,a_2\in G_{1_{G}},$ if $\Phi_{1_G}(a_1)=\Phi_{1_{G}}(a_2)$, then
$$a_1\circ 1_{G}=\mathcal{B}_1(a_1)\mathcal{B}_2(a_1)=\Phi_{1_G}(a_1)=\Phi_{1_{G}}(a_2)=\mathcal{B}_1(a_2)\mathcal{B}_2(a_2)=a_2\circ 1_G.$$
Hence
$$a_1=a_1\circ(1_G\circ 1^{\dag}_G)=(a_1\circ 1_G)\circ 1^{\dag}_G=(a_2\circ 1_G)\circ 1^{\dag}_G=a_2,$$
which implies that $\Phi_{1_G}$ is injective.
\qed

As a consequence of the above proposition, we have

\begin{cor}\label{B12}
Let $(G,\mathcal{B}_1,\mathcal{B}_2)$ be a Rota-Baxter system. Then for any $a\in G,$ $\Phi(a)=1_G$ if and only if $\mathcal{B}_1(a)=1_G$ and $\mathcal{B}_2(a)=1_G.$
\end{cor}

\proof
Assume that $\Phi(a)=1_G$, then by Theorem \ref{SEM},
$$\Phi(a\circ e_{1_G})=B_1(a\circ e_{1_G})B_2(a\circ e_{1_G})=B_1(a)B_2(a)=\Phi(a)=1_G=\Phi(e_{1_G}).$$
By Proposition \ref{IMG}, $\Phi_{1_G}$ is bijective. Then we have $a\circ e_{1_G}=e_{1_G}$. Hence we obtain that
$$\mathcal{B}_1(a)=\mathcal{B}_1(a\circ e_{1_G})=1_G,\quad \mathcal{B}_2(a)=\mathcal{B}_2(a\circ e_{1_G})=1_G,$$
as required.
\qed

Now we prove that $G$ is the disjoint union of the groups $G_t$, $t\in G$.

\begin{prop}\label{STR}
Let $(G,\mathcal{B}_1,\mathcal{B}_2)$ be a Rota-Baxter system. Then we have $$G=\bigsqcup {G_t}.$$
\end{prop}

\proof
For any $t_1,t_2\in G$, if there is $a\in G_{t_1}\cap G_{t_2}$, then we have
$$a\circ e_{t_1}=a=a\circ e_{t_2}.$$
Hence by Proposition \ref{INV}, we get $e_{t_1}=e_{t_2}$. Therefore using Lemma \ref{GT}, we find
$$G_{t_1}=G\circ e_{t_1}=G\circ e_{t_2}=G_{t_2},$$
which implies that $\bigcup G_t$ is a disjoint union.

As $t\in G_t$ for any $t\in G$, we have $G=\bigsqcup {G_t}$.
\qed

Set
$$K=\{e_{t}\mid t\in G\},$$
which is a subsemigroup of $(G,\circ)$. In fact, for any $a,b\in G$, we have $e_a\circ e_b=e_b$ by Theorem \ref{SEM}. The following theorem is the main result of this section.

\begin{thm}\label{OPL}
Let $(G,\mathcal{B}_1,\mathcal{B}_2)$ be a Rota-Baxter system. Then as semigroups,
$$(G,\circ)\cong (G_{1_G}\oplus K,\circ),$$
where the multiplication of $G_{1_G}\oplus K$ is given by
$$(a_1,b_1)\circ (a_2,b_2)=(a_1\circ a_2,b_1\circ b_2),\quad\forall a_1,a_2\in G_{1_G},\quad\forall b_1,b_2\in K.$$
\end{thm}

\proof
Define $\Psi:G\to G_{1_G}\oplus K$ by $\Psi(a)=(a\circ e_{1_G},e_a).$ We show $\Psi$ is an isomorphism of semigroups.

For any $a,b\in G,$ we have
$$\Psi(a)\circ \Psi(b)=(a\circ e_{1_G}\circ b\circ e_{1_G},e_a\circ e_b)
=(a\circ b\circ e_{1_G},e_b)=\Psi(a\circ b),$$
which proves that $\Psi$ is a homomorphism of semigroups.

For any $a\in G_{1_G}$ and $b\in G$, set $c=a\circ e_b$. By Lemma \ref{GT}, $c\in G_b$. Then we get
$$c\circ e_b=c=c\circ e_c,$$
which together with Proposition \ref{INV} implies that $e_b=e_c$. As
$$c\circ e_{1_G}=a\circ e_b\circ e_{1_G}=a\circ e_{1_G}=a,$$
we have $\Psi(c)=(a,e_b)$. Hence $\Psi$ is surjective.

It remains to show that $\Psi$ is injective. For any $a,b\in G,$ if $\Psi(a)=\Psi(b),$ then $a\circ e_{1_G}=b\circ e_{1_G}$ and $e_a=e_b$. Therefore
$$a=a\circ e_a=a\circ e_b=a\circ(e_{1_G}\circ e_b)=(a\circ e_{1_G})\circ e_b=(b\circ e_{1_G})\circ e_b=b,$$
which implies $\Psi$ is injective.
\qed

Inspired by Theorem \ref{OPL}, there is a way to construct Rota-Baxter systems.

\begin{eg}
Let $(G,\mathcal{B})$ be a Rota-Baxter group. Define a map $\widetilde{\mathcal{B}}:G\to G$ by
$$\widetilde{\mathcal{B}}(a)=a\mathcal{B}(a),\quad\forall a\in G.$$
It follows from \cite[Remark 2.11, Proposition 3.1]{LG1} that $\widetilde{\mathcal{B}}$ is a Rota-Baxter operator of weight $-1$. Let $K$ be another group. Define $\mathcal{\mathcal{B}}_1:G\oplus K\to G\oplus K$ by $\mathcal{B}_1((a,b))=(\widetilde{\mathcal{B}}(a),1_K)$ and $\mathcal{B}_2:G\oplus K\to G\oplus K$ by $\mathcal{B}_2((a,b))=(\mathcal{B}(a)^{-1},1_K)$ for any $a\in G$ and $b\in K$. One can verify that $(G\oplus K,\mathcal{B}_1,\mathcal{B}_2)$ is a Rota-Baxter system.
\end{eg}

For a Rota-Baxter system $(G,\mathcal{B}_1,\mathcal{B}_2)$, define
\begin{equation}
\Psi:G\to G_{1_G},\quad \Psi(a)=a\circ e_{1_G},\quad \forall a\in G. \label{Psi1}
 \end{equation}
 Then for any $a\in G$, we have
$$\mathcal{B}_1(\Psi(a))=\mathcal{B}_1(a\circ e_{1_G})=\mathcal{B}_1(a)\mathcal{B}_1(e_{1_G})=\mathcal{B}_1(a)$$
and
$$\mathcal{B}_2(\Psi(a))=\mathcal{B}_2(a\circ e_{1_G})=\mathcal{B}_2(e_{1_G})\mathcal{B}_2(a)=\mathcal{B}_2(a).$$
Furthermore, one can verify  $\Phi_{1_G}^{-1}(\Phi(a))=\Psi(a).$

\begin{prop}\label{ROT}
Let $(G,\mathcal{B}_1,\mathcal{B}_2)$ be a Rota-Baxter system. Then there are two operators $\mathcal{C}_1:G\to G$ and $\mathcal{C}_2:G\to G$ such that $\mathcal{C}_1(1_G)=\mathcal{C}_2(1_G)=1_G$ and $(G,\mathcal{C}_1,\mathcal{C}_2)$ is a Rota-Baxter system.
\end{prop}

\proof
We define
$$\mathcal{C}_1(a)=\mathcal{B}_1(\Phi_{1_G}^{-1}(\Psi(a))),\quad \mathcal{C}_2(a)=\mathcal{B}_2(\Phi_{1_G}^{-1}(\Psi(a))),\quad \forall a\in G.$$
For any $a,b\in G,$ we have
$$
\mathcal{C}_1(a)\mathcal{C}_1(b)=\mathcal{B}_1(\Phi_{1_G}^{-1}(\Psi(a))\circ \Phi_{1_G}^{-1}(\Psi(b)))
=\mathcal{B}_1\Phi_{1_G}^{-1}\Phi(\Phi_{1_G}^{-1}(\Psi(a))\circ \Phi_{1_G}^{-1}(\Psi(b))).$$
Using Lemma \ref{TB},
\begin{align*}
\mathcal{C}_1(a)\mathcal{C}_1(b)&=\mathcal{B}_1\Phi_{1_G}^{-1}(\Phi_{1_G}^{-1}(\Psi(a))\circ \Psi(b))
=\mathcal{B}_1\Phi_{1_G}^{-1}(\Phi_{1_G}^{-1}(\Psi(a))\circ b\circ e_{1_G})\\
&=\mathcal{B}_1\Phi_{1_G}^{-1}\Psi(\Phi_{1_G}^{-1}(\Psi(a))\circ b)=\mathcal{C}_1(\mathcal{B}_1(\Phi_{1_G}^{-1}(\Psi(a))) b\mathcal{B}_2(\Phi_{1_G}^{-1}(\Psi(a))))\\
&=\mathcal{C}_1(\mathcal{C}_1(a)b\mathcal{C}_2(b)).
\end{align*}
One can verify $\mathcal{C}_2(b)\mathcal{C}_2(a)=\mathcal{C}_2(\mathcal{C}_1(a)b\mathcal{C}_2(a))$ similarly. Hence $(G,\mathcal{C}_1,\mathcal{C}_2)$ is a Rota-Baxter system.

We have
$$\mathcal{C}_1(1_G)=\mathcal{B}_1(\Phi_{1_G}^{-1}(\Psi(1_G)))=\mathcal{B}_1(\Phi_{1_G}^{-1}(1_G))=\mathcal{B}_1(e_{1_G})=1_G.$$
Similarly, we can show $\mathcal{C}_2(1_G)=1_G$.
\qed

We finally show that the operator $\Psi$ is the cocycle of the Rota-Baxter system $(G,\mathcal{C}_1,\mathcal{C}_2)$ in Proposition \ref{ROT}.

\begin{lem}\label{LEM}
Let $(G,\mathcal{B}_1,\mathcal{B}_2)$ be a Rota-Baxter system such that $\mathcal{B}_1(1_G)=\mathcal{B}_2(1_G)=1_G.$ Then the cocycle $\Phi$ of $(G,\mathcal{B}_1,\mathcal{B}_2)$ is just the operator $\Psi:G\to G$ defined by \eqref{Psi1}. Moreover
\begin{equation}
\mathcal{B}_1(\Phi(a))=\mathcal{B}_1(a),\quad \mathcal{B}_2(\Phi(a))=\mathcal{B}_2(a)
\end{equation} for any $a\in G.$
\end{lem}

\proof
For any $a\in G,$
$$e_{1_G}=\mathcal{B}_1(1_G)^{-1}1_G\mathcal{B}_2(1_G)^{-1}=1_G.$$
Hence
$$\Psi(a)=a\circ e_{1_G}=a\circ 1_G=\mathcal{B}_1(a)\mathcal{B}_2(a)=\Phi(a),$$
which implies $\Phi=\Psi$.
\qed

%\section{Right quasi-invertible semigroup and embedding of skew truss}
%This section is divided into two parts. In the first part, the notion of right quasi-invertible semigroup is introduced and some results in \cite{VG} are generalized. In the second, the embedding from skew truss to Rota-Baxter system is constructed.
%\subsection{Right quasi-invertible semigroup and $\Phi$-holomorph}
%\begin{Def}
%Let $(G,\circ)$ be a semigroup. $(G,\circ)$ is called right quasi-invertible or quasi-invertible in short if for any $a,b\in G$, there exists a $c\in G$ such that $a\circ c=b$.
%\end{Def}
%From Proposition \ref{INV}, if $\circ$ is an associative descendent operation of $(G,\mathcal{B}_1,\mathcal{B}_2)$, then the semigroup $(G,\circ)$ is quasi-invertible.

 \section{Factorization of Rota-Baxter systems}\label{5}

In this section, we prove a factorization theorem of Rota-Baxter systems, which is a generalization of \cite[Theorem 3.5]{LG1}. The theorem can be regarded as an analogue of the global factorization theorem of M. A. Semenov-Tian-Shansky for Lie groups given in \cite{Semenov}.

Let $(G,\mathcal{B}_1,\mathcal{B}_2)$ be a Rota-Baxter system with the descendent operation $\circ$ and the cocycle $\Phi$. As before, for any $a\in G$, set $e_a=\mathcal{B}_1(a)^{-1}a\mathcal{B}_2(a)^{-1}$ and $a^{\dag}=\mathcal{B}_1(a)^{-1}e_a\mathcal{B}_2(a)^{-1}$. Then by Theorem \ref{SEM}, $e_a$ is a left identity with respect to $\circ$, and by Lemma \ref{DAG}, $\mathcal{B}_1(a^{\dag})=\mathcal{B}_1(a)^{-1}$ and $\mathcal{B}_2(a^{\dag})=\mathcal{B}_2(a)^{-1}$. Moreover, we define four subsets of $G$:
$$G_1=\operatorname{Im}(\mathcal{B}_1),\quad G_2=\operatorname{Im}(\mathcal{B}_2), \quad H_1=\mathcal{B}_1(\operatorname{Ker}(\mathcal{B}_2)),\quad H_2=\mathcal{B}_2(\operatorname{Ker}(\mathcal{B}_1)).$$
Then $G_1$ and $G_2$ are subgroups of $G$ by Corollary \ref{INV1}.

\begin{lem}
$H_1$ is a subgroup of $G_1$ and $H_2$ is a subgroup of $G_2$.
\end{lem}

\proof
For any $a,b\in \operatorname{Ker}(\mathcal{B}_1)$, we have
$$\mathcal{B}_1(b\circ a)=\mathcal{B}_1(b)\mathcal{B}_1(a)=1_G,$$
which induces that
$$\mathcal{B}_2(a)\mathcal{B}_2(b)=\mathcal{B}_2(b\circ a)\in H_2.$$
Hence $H_2$ is closed under the multiplication of $G_2.$ Using Theorem \ref{SEM}, we have
$$\mathcal{B}_1(e_a)=\mathcal{B}_2(e_a)=1_G,$$
which implies $1_G\in H_2$. As
$$\mathcal{B}_1(a)\mathcal{B}_1(a^{\dag})=1_G,$$
we find $a^{\dag}\in \operatorname{Ker}(\mathcal{B}_1)$. Hence $\mathcal{B}_{2}(a)^{-1}=\mathcal{B}_2(a^{\dag})\in H_2$. Then $H_2$ is a subgroup of $G_2$.

Similarly, we can prove that $H_1$ is a subgroup of $G_1$.
\qed

In fact, $H_i$ is a normal subgroup of $G_i$. To prove this result, we need a lemma.

\begin{lem}\label{INV2}
 Let $(G,\mathcal{B}_1,\mathcal{B}_2)$ be a Rota-Baxter system. Then
\begin{description}
  \item [(1)]  for any $a\in G$ and $b\in \operatorname{Ker}(\mathcal{B}_1)$, there exists $c\in \operatorname{Ker}(\mathcal{B}_1)$ such that $\mathcal{B}_2(c)\mathcal{B}_2(a)=\mathcal{B}_2(a)\mathcal{B}_2(b)$;
  \item [(2)] for any $a_1\in G$ and $b_1\in \operatorname{Ker}(\mathcal{B}_2)$, there exists $c_1\in \operatorname{Ker}(\mathcal{B}_2)$ such that $\mathcal{B}_1(a_1)\mathcal{B}_1(c_1)=\mathcal{B}_1(b_1)\mathcal{B}_1(a_1)$.
\end{description}
\end{lem}

\proof
(1) By the proof of Proposition \ref{INV}, taking $c=\mathcal{B}_1(a)^{-1}(b\circ a)\mathcal{B}_2(a)^{-1}$, we have $a\circ c=b\circ a$. As $b\in \operatorname{Ker}(\mathcal{B}_1)$, we get
$$\mathcal{B}_1(a)\mathcal{B}_1(c)=\mathcal{B}_1(a\circ c)=\mathcal{B}_1(b\circ a)=\mathcal{B}_1(a).$$
Then $\mathcal{B}_1(c)=1_G$ and $c\in \operatorname{Ker}(\mathcal{B}_1)$. Moreover, we have
$$\mathcal{B}_2(c)\mathcal{B}_2(a)=\mathcal{B}_2(a\circ c)=\mathcal{B}_2(b\circ a)=\mathcal{B}_2(a)\mathcal{B}_2(b).$$

(2) Similarly, taking $c_1=\mathcal{B}_1(a)^{-1}(b_1\circ a_1)\mathcal{B}_2(a)^{-1},$ we obtain $a_1\circ c_1=b_1\circ a_1$. As $b_1\in \operatorname{Ker}(\mathcal{B}_2)$ and $$\mathcal{B}_2(c_1)\mathcal{B}_2(a_1)=\mathcal{B}_2(a_1\circ c_1)=\mathcal{B}_2(b_1\circ a_1)=\mathcal{B}_2(a_1),$$
we have $c_1\in \operatorname{Ker}(\mathcal{B}_2)$. Finally, we derive that
$$\mathcal{B}_1(a_1)\mathcal{B}_1(c_1)=\mathcal{B}_1(a_1\circ c_1)=\mathcal{B}_1(b_1\circ a_1)=\mathcal{B}_1(b_1)\mathcal{B}_1(a_1),$$
which finishes the proof.
\qed

\begin{lem}\label{NORM}
$H_1$ is a normal subgroup of $G_{1}$ and $H_2$ is a normal subgroup of $G_{2}.$
\end{lem}

\proof
Using Lemma \ref{INV2} (1), for any $a\in G$, we have $\mathcal{B}_2(a)H_2\subseteq H_2\mathcal{B}_2(a)$. Hence $H_2$ is a normal subgroup of $G_2$. Similarly, Lemma \ref{INV2} (2) implies that $H_1$ is a normal subgroup of $G_1$. \qed

Based on Lemma \ref{NORM}, we define the map:
$$\Theta:G_1/H_1\to G_2/H_2,\quad \Theta(\overline{\mathcal{B}_1(a)})=\overline{\mathcal{B}_2(a)},\quad \forall a\in G.$$
where $\overline{\cdot}$ denotes the equivalence classes in the two quotients. As for any $b\in \operatorname{Ker}(\mathcal{B}_2)$,
$$\Theta(\overline{\mathcal{B}_1(a)\mathcal{B}_1(b)})=\Theta(\overline{\mathcal{B}_1(a\circ b)})=\overline{\mathcal{B}_2(a\circ b)}=\overline{\mathcal{B}_2(b)\mathcal{B}_2(a)}
=\overline{\mathcal{B}_2(a)}=\Theta(\overline{\mathcal{B}_1(a)}),$$
the map $\Theta$ is well-defined.

\begin{prop}\label{CAL}
The map $\Theta:G_1/H_1\to G_2/H_2$ is a group anti-isomorphism. It is called the Cayley transform of the Rota-Baxter system $(G,\mathcal{B}_1,\mathcal{B}_2).$
\end{prop}

\proof
For any $a, b\in G,$  we have
$$\Theta(\overline{\mathcal{B}_1(a)\mathcal{B}_1(b)})=\Theta(\overline{\mathcal{B}_1(a\circ b)})=\overline{\mathcal{B}_2(a\circ b)}=\overline{\mathcal{B}_2(b)\mathcal{B}_2(a)}=\Theta(\overline{\mathcal{B}_1(b)})\Theta(\overline{\mathcal{B}_1(a)}),$$
which implies $\Theta$ is a group anti-homomorphism.

It's obvious that $\Theta$ is surjective. Hence it remains to show that $\Theta$ is injective. If $\mathcal{B}_2(a)=\mathcal{B}_2(b),$ where $a\in G$ and $b\in \operatorname{Ker}(\mathcal{B}_1)$, then we have
$$\mathcal{B}_2(b^{\dag}\circ a)=\mathcal{B}_2(a)\mathcal{B}_2(b^{\dag})=\mathcal{B}_2(a)\mathcal{B}_2(b)^{-1}=1_G,$$
which implies that $b^{\dag}\circ a\in \operatorname{Ker}(\mathcal{B}_2)$. We next derive that
$$\mathcal{B}_1(a)=\mathcal{B}_1(b^{\dag})\mathcal{B}_1(a)=\mathcal{B}_1(b^{\dag}\circ a),$$
which proves $\Theta$ is injective.
\qed

Consider the group $(G_1\times G_2,\cdot_D)$, where the product is defined by
$$(a_1,a_2)\cdot_D (b_1,b_2)=(a_1b_1,b_2a_2),\quad \forall a_1,b_1\in G_1,a_2,b_2\in G_2. $$
Let $G_\Theta$ denote the subset
$$G_{\Theta}=\{(a_1,a_2)\in G_1\times G_2\mid \Theta(\overline{a_1})=\overline{a_2}\}.$$
Define a map $\Psi:G_{1_G}\to G_{\Theta}$ by
\begin{equation}\label{Psi}
\Psi(a)=(\mathcal{B}_1(a),\mathcal{B}_2(a)),\quad \forall a\in G_{1_G}.
\end{equation}

\begin{lem}\label{G1}
$G_{\Theta}$ is a subgroup of $(G_1\times G_2,\cdot_D)$ and $\Psi$ is a group isomorphism from $(G_{1_G},\circ)$ to $(G_1\times G_2,\cdot_D)$.
\end{lem}

\proof
Using Proposition \ref{CAL}, for any $(a_1,a_2),(b_1,b_2)\in G_{\Theta}$, we have $$\Theta(\overline{a_1b_1})=\Theta(\overline{a_1}\overline{b_1})=\Theta(\overline{b_1})\Theta(\overline{a_1})=\overline{b_2}\overline{a_2}=\overline{b_2a_2},$$
which implies that $(a_1b_1,b_2a_2)\in G_{\Theta}$. Similarly, it is easy to see that $(a_1,b_1)^{-1}\in G_{\Theta}$. Hence $G_{\Theta}$ is a subgroup of $G_1\times G_2.$

For any $a,b\in G_{1_G}$,
\begin{align*}
\Psi(a\circ b)&=(\mathcal{B}_1(a\circ b),\mathcal{B}_2(a\circ b))=(\mathcal{B}_1(a)\mathcal{B}_1(b),\mathcal{B}_2(b)\mathcal{B}_2(a))\\
&=(\mathcal{B}_1(a),\mathcal{B}_2(a))\cdot_D (\mathcal{B}_1(b),\mathcal{B}_2(b))=\Psi(a)\cdot_D \Psi(b),
\end{align*}
which implies that $\Psi$ is a homomorphism.

For an element $a\in G_{1_G}$ such that $\Psi(a)=(1_G,1_G)$, by Corollary \ref{B12},
$$\Phi_{1_G}(a)=1_G=\Phi_{1_G}(e_{1_G}).$$
Hence $a=e_{1_G}$, and $\Psi$ is injective.

For any $(a_1,a_2)\in G_{\Theta}$, we have $\Theta(\overline{a_1})=\overline{a_2}$. There exists  $a\in G$ such that $a_1=\mathcal{B}_1(a)$. Hence $$\Theta(\overline{a_1})=\Theta(\overline{\mathcal{B}_1(a)})=\overline{\mathcal{B}_2(a)},$$
which implies that $\overline{a_2}=\overline{\mathcal{B}_2(a)}$. Therefore there exists $b\in \operatorname{Ker}(\mathcal{B}_1)$ such that
$$a_2=\mathcal{B}_2(a)\mathcal{B}_2(b).$$
Set $a_3=b\circ a\circ e_{1_G}\in G_{1_G}$, then
$$\Psi(a_3)=(\mathcal{B}_1(b\circ a),\mathcal{B}_2(b\circ a))=(\mathcal{B}_1(b)\mathcal{B}_1(a),\mathcal{B}_2(a)\mathcal{B}_2(b))=(a_1,a_2),$$
which means that $\Psi$ is surjective.
\qed

Finally, we give the following decomposition theorem of Rota-Baxter systems.

\begin{thm}\label{Dec}
Let $(G,\mathcal{B}_1,\mathcal{B}_2)$ be a Rota-Baxter system. Then every element $a\in \Phi(G)=G_{1_G}$ can be uniquely decomposed as $a=a_1a_2$ with $(a_1,a_2)\in G_{\Theta}$.
\end{thm}

\proof
If $a=\Phi(b)\in G_{1_G}$ with $b\in G$, set $a_1=\mathcal{B}_1(b)$ and $a_2=\mathcal{B}_2(b)$. Then we have $(a_1,a_2)\in G_{\Theta}$ and
$$a=\Phi(b)=a_1a_2.$$
It remains to show the uniqueness of the decomposition. If $a=a_1a_2=b_1b_2$ with $(a_1,a_2),(b_1,b_2)\in G_{\Theta}$, then we have $a_1^{-1}b_1b_2a_2^{-1}=1_G$. Since $(G_{\Theta},\cdot_D)$ is a group, we have $(a_1^{-1}b_1,b_2a_2^{-1})\in G_{\Theta}$. By Lemma \ref{G1}, there exists $t\in G_{1_G}$ such that
$$\Psi(t)=(\mathcal{B}_1(t),\mathcal{B}_2(t))=(a_1^{-1}b_1,b_2a_2^{-1}),$$
where $\Psi$ is defined by \eqref{Psi}. Therefore
$$\Phi(t)=\mathcal{B}_1(t)\mathcal{B}_2(t)=a_1^{-1}b_1b_2a_2^{-1}=1_G.$$
Using Proposition \ref{B12}, we have $\mathcal{B}_1(t)=\mathcal{B}_2(t)=1_G$. Hence $a_1=b_1, a_2=b_2$.
 \qed

\begin{rem}
Theorem \ref{Dec} describes the decomposition of $G_{1_G}$, which is different from the global decomposition theorem in \cite[Theorem 3.5]{LG1}. While from Lemma \ref{ISO} and Proposition \ref{STR}, under a group isomorphism, every element of $G$ can be mapped into $G_{1_G}$, and then be decomposed by Theorem \ref{Dec}.
\end{rem}

\section{Rota-Baxter systems of Lie algebras and twisted modified Yang-Baxter equations}\label{6}

In this section, we study the Rota-Baxter systems of Lie algebras. We first propose the notion of Rota-Baxter systems of Lie algebras and show their factorization theorem. Then we study the relationship between Rota-Baxter systems of Lie algebras and of Lie groups. We show that Rota-Baxter systems of Lie groups are the integration of Rota-Baxter systems of Lie algebras under some conditions. Finally, the notion of twisted modified Yang-Baxter equations is introduced as a generalization of the modified Yang-Baxter equation. It is shown that the twisted modified Yang-Baxter equations have a solution on a Lie algebra if and only if the Lie algebra has a structure of Rota-Baxter systems.

\subsection{Rota-Baxter systems of Lie algebras}

In this subsection, we introduce the notion of Rota-Baxter systems of Lie algebras. And we prove a factorization theorem as a generalization of the global factorization theorem in \cite{Semenov}.

\begin{Def}
A triple $(\mathfrak{g},B_1,B_2)$ with a Lie algebra $\mathfrak{g}$ and two linear operators $B_1:\mathfrak{g}\to \mathfrak{g}$ and $B_2:\mathfrak{g}\to \mathfrak{g}$ is called a Rota-Baxter system of Lie algebras if for any $u,v\in \mathfrak{g}$,
\begin{subequations}
\begin{align}
[B_1(u),B_1(v)]&=B_1([B_1(u),B_1(v)]-[B_2(u),B_2(v)]),\label{TWT1}\\
[B_2(v),B_2(u)]&=B_2([B_1(u),B_1(v)]-[B_2(u),B_2(v)]).\label{TWT2}
\end{align}
\end{subequations}
\end{Def}

Taking $\phi=B_1+B_2,$ we may rewrite \eqref{TWT1} and \eqref{TWT2} as
\begin{subequations}
\begin{align}
[B_1(u),B_1(v)]&=B_1([B_1(u),\phi(v)]+[\phi(u),B_1(v)]-[\phi(u),\phi(v)]), \label{RBSL1}\\
[B_2(u),B_2(v)]&=B_2([B_2(u),\phi(v)]+[\phi(u),B_2(v)]-[\phi(u),\phi(v)]).
\end{align}
\end{subequations}
Moreover, if $\phi=\operatorname{id}$, then $B_1$ and $B_2$ are Rota-Baxter operators of Lie algebras of weight $-1$. Hence Rota-Baxter systems of  Lie algebras generalize the notion of Rota-Baxter Lie algebras.

The next example and proposition generalize \cite[Proposition 5]{Semenov}.

\begin{eg}\label{EgRBS}
Let $\mathfrak{g}$ be a Lie algebra with subalgebras $\mathfrak{g}_\pm$ and a subspace $V$.  Assume that $\mathfrak{g}=\mathfrak{g}_+\oplus\mathfrak{g}_-\oplus V$ as vector spaces. Let $B_{1}$ be the projection from $\mathfrak{g}$ to $\mathfrak{g}_{+}$ and $B_2$ be the projection from $\mathfrak{g}$ to $\mathfrak{g}_{-}$. Then it is easy to see that $(\mathfrak{g},B_1,B_2)$ is a Rota-Baxter system of Lie algebras.
\end{eg}

\begin{prop}\label{EOR}
Let $(\mathfrak{g},B_1,B_2)$ be Rota-Baxter system of Lie algebras. Then the operation
\begin{equation}
[\cdot,\cdot]_R:\mathfrak{g}\times\mathfrak{g}\to \mathfrak{g},\quad [u,v]_R=[B_1(u),B_1(v)]-[B_2(u),B_2(v)]\label{LBR},
\end{equation}
is a Lie bracket. Denote by $\mathfrak{g}_R$ the Lie algebra $\mathfrak{g}$ with respect to the Lie bracket $[\cdot,\cdot]_R$.
\end{prop}

\proof
It's enough to show $[\cdot,\cdot]_R$ satisfies the Jacobi identity. For any $u,v,w\in\mathfrak{g}$, by \eqref{TWT1} and \eqref{TWT2}, we have
\begin{align*}
[u,[v,w]_R]_R=&[B_1(u),B_1([v,w]_R)]-[B_2(u),B_2([v,w]_R)]\\
=&[B_1(u),[B_1(v),B_1(w)]]+[B_2(u),[B_2(v),B_2(w)]].
\end{align*}
Hence we get
$$[u,[v,w]_R]_R+[v,[w,u]_R]_R+[w,[u,v]_R]_R=0,$$
as required.
\qed

\begin{prop}
Let $(\mathfrak{g},B_1,B_2)$ be a Rota-Baxter system of Lie algebras. Then
\begin{description}
 \item[(1)] $B_1:\mathfrak{g}_R\to \mathfrak{g}$ is a Lie algebra homomorphism and $B_2:\mathfrak{g}_R\to \mathfrak{g}$ is a Lie algebra anti-homomorphism;
 \item[(2)] $B_1(\mathfrak{g})\subseteq \mathfrak{g}$ and $B_2(\mathfrak{g})\subseteq\mathfrak{g}$ are Lie subalgebras;
 \item[(3)] $\operatorname{Ker}(B_1)\subseteq \mathfrak{g}_R$ and $\operatorname{Ker}(B_2)\subseteq \mathfrak{g}_R$ are ideals, and there are Lie algebra isomorphism $B_1(\mathfrak{g})\cong\mathfrak{g}_R/\operatorname{Ker}(B_1)$ and Lie algebra anti-isomorphism $B_2(\mathfrak{g})\cong\mathfrak{g}_R/\operatorname{Ker}(B_2)$.
\end{description}
\end{prop}

\proof
Using \eqref{TWT1} and \eqref{TWT2}, we get (1). Then (2) and (3) are straightforward consequences of (1).
\qed

As a generalization of \cite[Proposition 8]{Semenov}, the Cayley transform $\theta$ is established in the following proposition.

\begin{prop}
Let $(\mathfrak{g},B_1,B_2)$ be a Rota-Baxter system of Lie algebras. Then
\begin{description}
  \item[(1)] $B_1(\operatorname{Ker}(B_2))\subseteq B_1(\mathfrak{g})$ and $B_2(\operatorname{Ker}(B_1))\subseteq B_2(\mathfrak{g})$ are ideals;
  \item[(2)] the map $\theta:B_1(\mathfrak{g})/B_1(\operatorname{Ker}(B_2))\to B_2(\mathfrak{g})/B_2(\operatorname{Ker}(B_1))$ defined by
  $$\theta(\overline{B_1(u)})=\overline{B_2(u)}$$
  is a Lie algebra anti-isomorphism.
\end{description}
\end{prop}

\proof
The statement (1) follows from \eqref{TWT1} and \eqref{TWT2}.

For any $u\in \mathfrak{g}$ and $v\in \operatorname{Ker}(B_2)$, we have
$$\overline{B_2(u+v)}=\overline{B_2(u)},$$
which implies $\theta$ is well-defined. For any $u,v\in\mathfrak{g}$, by \eqref{TWT1} and \eqref{TWT2},
\begin{align*}
\theta(\overline{[B_1(u),B_1(v)]})=\theta(\overline{B_1([u,v]_R)})
=\overline{B_2([u,v]_R)}=\overline{[B_2(v),B_2(u)]},
\end{align*}
which proves $\theta$ is a Lie algebra anti-homomorphism. It is obvious that $\theta$ is surjective. Hence it remains to show it's injective. If $u\in \mathfrak{g}$ and $v\in \operatorname{Ker}(B_1)$ such that $B_2(u)=B_2(v)$, then $u-v\in \operatorname{Ker}(B_2)$. Hence
$$B_1(u)=B_1(u-v)+B_1(v)=B_1(u-v)\in B_1(\operatorname{Ker}(B_2)),$$
which proves $\theta$ is injective. Therefore we get (2).
\qed

For a Rota-Baxter system of Lie algebras $(\mathfrak{g},B_1,B_2)$, set $\phi=B_1+B_2$.

\begin{lem}
$\phi(\mathfrak{g})\subseteq \mathfrak{g}_R$ is a Lie subalgebra.
\end{lem}

\proof
For any $u,v\in \mathfrak{g},$ we have $$\phi([u,v]_R)=B_1([u,v]_R)+B_2([u,v]_R)=[B_1(u),B_1(v)]-[B_2(u),B_2(v)]=[u,v]_R.$$
Hence $\phi(\mathfrak{g})$ is closed under $[\cdot,\cdot]_R$.
\qed

Consider
$$\tilde{\mathfrak{g}}=\{(u_+,u_-)\in B_1(\mathfrak{g})\oplus B_2(\mathfrak{g})\mid \theta(\overline{u_+})=\overline{u_-}\},$$
which is a subalgebra of $B_1(\mathfrak{g})\oplus B_2(\mathfrak{g})$. Define a map
$$T:\phi(\mathfrak{g})\to \tilde{\mathfrak{g}},\quad T(u)=(B_1(u),B_2(u)).$$
The following proposition generalizes \cite[Proposition 9]{Semenov}.

\begin{prop}
Let $(\mathfrak{g},B_1,B_2)$ be a Rota-Baxter system of Lie algebras. If $\phi$ is idempotent, then $T$ is a Lie algebra isomorphism.
\end{prop}

\proof
For any $u,v\in \mathfrak{g}$, we have
\begin{align*}
[T(\phi(u)),T(\phi(v))]&=\left([B_1(\phi(u)),B_1(\phi(v))], [B_2(\phi(u)),B_2(\phi(v))]\right)\\
&=\left(B_1([\phi(u),\phi(v)]_R), B_2([\phi(u),\phi(v)]_R)\right)\\
&=T([\phi(u),\phi(v)]_R),
\end{align*}
which implies $T$ is a Lie algebra homomorphism.

Let $u\in \mathfrak{g}$ such that $T(\phi(u))=0,$ then $B_1(\phi(u))=B_2(\phi(u))=0$. Hence
$$\phi(u)=\phi^{2}(u)=B_1(\phi(u))+B_2(\phi(u))=0,$$
which implies $T$ is injective.

For any $(u_+,u_-)\in \tilde{\mathfrak{g}}$, we have $\theta(\overline{u_+})=\overline{u_-}$. Since $u_+\in B_1(\mathfrak{g}),$ there is $u\in \mathfrak{g}$ such that $B_1(u)=u_+.$ Therefore we obtain
$$\theta(\overline{u_+})=\theta(\overline{B_1(u)})=\overline{B_2(u)}.$$
Hence $\overline{u_-}=\overline{B_2(u)}$. Then there exists $v\in \operatorname{Ker}(B_1)$ such that $B_2(u)=u_- + B_2(v).$ Taking $u_1=u -v,$ we have
$$T(u_1)=(B_1(u_1),B_2(u_1))=(B_1(u)-B_1(v),B_2(u)-B_2(v))=(u_+,u_-).$$
Therefore, $T$ is surjective.
\qed

Now we give a factorization theorem of Rota-Baxter systems of Lie algebras.

\begin{Def}
A Rota-Baxter system of Lie algebras $(\mathfrak{g},B_1,B_2)$ is called decomposable if $\operatorname{Ker}(\phi)\subseteq \operatorname{Ker}(B_1).$
\end{Def}

In the above definition, we may replace $\operatorname{Ker}(B_1)$ by $\operatorname{Ker}(B_2)$ or $\operatorname{Ker}(B_1)\cap \operatorname{Ker}(B_2)$.

\begin{thm}[Factorization theorem of Rota-Baxter systems of Lie algebras]
Let $(\mathfrak{g},B_1,B_2)$ be a decomposable Rota-Baxter system of Lie algebras. Then any $\phi(u)\in\phi(\mathfrak{g})$ can be uniquely decomposed as $\phi(u)=u_++u_-$ with $(u_+,u_-)\in \tilde{\mathfrak{g}}.$
\end{thm}

\proof
For any $u\in \mathfrak{g},$ $\phi(u)=B_1(u)+B_2(u).$ As $(B_1(u),B_2(u))\in\tilde{\mathfrak{g}}$, we get the existence of the decomposition. To see the uniqueness of the decomposition, let $\phi(u)=u_++u_-=v_++v_-$ with $(u_+,u_-),(v_+,v_-)\in\tilde{\mathfrak{g}}$. Then $u_+-v_+= v_--u_-\in B_1(\mathfrak{g})\cap B_2(\mathfrak{g})$ and $\theta(\overline{u_+-v_+})=\overline{u_- -v_-}$. There exists $w\in \mathfrak{g}$ such that
$$u_+-v_+= v_- - u_-=B_1(w)\in B_1(\mathfrak{g})\cap B_2(\mathfrak{g}).$$
Then
$$\theta(\overline{B_1(w)})=\overline{B_2(w)}=\overline{-B_1(w)}.$$
Hence $B_2(w)+B_2(w_1)=-B_1(w)$ for some $w_1\in \operatorname{Ker}(B_1)$, which implies $\phi(w+w_1)=0$. Since $\operatorname{Ker}(\phi)\subseteq \operatorname{Ker}{B_1}$,  we have
$$u_+-v_+=B_1(w)=B_1(w+w_1)=0.$$
Therefore $u_+=v_+$ and $u_-=v_-.$
\qed

One can readily verify that the Rota-Baxter system of Lie algebras defined in Example $\ref{EgRBS}$ is decomposable, hence it gives a class of examples of decomposable Rota-Baxter systems of Lie algebras.
In the following subsection, we will see that there is a closed connection between decomposable Rota-Baxter systems of Lie algebras and Rota-Baxter systems of groups, which is just like the relationship between Rota-Baxter Lie algebras and Rota-Baxter groups.

Finally, we give a lemma which will be used in the next subsection.

\begin{lem}\label{LIN}
Let $(\mathfrak{g},B_1,B_2)$ be a Rota-Baxter system of Lie algebras with $\phi=B_1+B_2$. Then the following statements are equivalent:
\begin{description}
 \item[(i)] $(\mathfrak{g},B_1,B_2)$ is decomposable and $\phi$ is idempotent;
 \item[(ii)] For any $u\in \mathfrak{g}$, $B_1(u)=B_1(\phi(u))$ and $B_2(u)=B_2(\phi(u))$.
\end{description}
\end{lem}

\proof
Suppose that $(\mathfrak{g},B_1,B_2)$ is decomposable and $\phi$ is idempotent. Then for any $u\in \mathfrak{g}$, we have $\phi(\phi(u)-u)=0$. Using $\operatorname{Ker}(\phi)\subseteq \operatorname{Ker}(B_1)\cap \operatorname{Ker}(B_2)$, we have
$$B_1(\phi(u)-u)=0,\quad B_2(\phi(u)-u)=0.$$
Therefore $B_1(u)=B_1(\phi(u))$ and $B_2(u)=B_2(\phi(u))$.

Conversely, suppose that for any $u\in \mathfrak{g}$, $B_1(u)=B_1(\phi(u))$ and $B_2(u)=B_2(\phi(u))$. Then
$$\phi(u)=B_1(u)+B_2(u)=\phi(\phi(u)),$$
which implies $\phi$ is idempotent. If $u\in \operatorname{Ker}(\phi)$, then
$$B_1(u)=B_1(0)=0,$$
which induces that $\operatorname{Ker}(\phi)\subseteq \operatorname{Ker}(B_1)$.
\qed

\subsection{Integrations of Rota-Baxter systems of Lie algebras}

We now show that Rota-Baxter systems of groups are integrations of decomposable Rota-Baxter systems of Lie algebras under some conditions.
Let $G$ be a lie group over $\mathbb{K}$. Let $\mathfrak{g}=T_{1_G}G$ be the Lie algebra of $G$ with Lie bracket $[\cdot,\cdot]$ and let
$$\exp:\mathfrak{g}\to G$$
be its exponential map. The Lie bracket $[\cdot,\cdot]$ and the group multiplication is connected by the formula:
$$[u,v]=\left.\frac{\mathrm{d}^{2}}{\mathrm{d}s \mathrm{d}t}\right|_{s,t=0}\exp(su)\exp(tv)\exp(-su),\quad \forall u,v\in \mathfrak{g}.$$

\begin{Def}
A Rota-Baxter system $(G,\mathcal{B}_1,\mathcal{B}_2)$ of groups is called a Rota-Baxter system of Lie groups, if it satisfies
\begin{description}
 \item[(1)] $G$ is a Lie group and $\mathcal{B}_1,\mathcal{B}_2$ are smooth operators on $G$;
 \item[(2)] $\mathcal{B}_1(1_G)$=$\mathcal{B}_2(1_G)=1_G.$
\end{description}
\end{Def}

By Proposition \ref{ROT}, if there is a Rota-Baxter system structure on $G,$ there always exist $\mathcal{C}_1:G\to G$ and $\mathcal{C}_2:G\to G,$ such that
$(G,\mathcal{C}_1,\mathcal{C}_2)$ is a Rota-Baxter system and $\mathcal{C}_1(1_G)$=$\mathcal{C}_2(1_G)=1_G$.

The next theorem generalizes \cite[Theorem 2.9]{LG1}. See also  \cite[Corrollary 2.14]{LG1}.

\begin{thm}\label{RBL}
Let $(G,\mathcal{B}_1,\mathcal{B}_2)$ be a Rota-Baxter system of Lie groups. Let $\mathfrak{g}$ be the Lie algebra of $G$ and $$B_1:\mathfrak{g}\to \mathfrak{g},\quad B_2:\mathfrak{g}\to \mathfrak{g}$$ be the
tangent map of $\mathcal{B}_1$ and $\mathcal{B}_2$ at $1_G$ respectively. Then for any $u,v\in\mathfrak{g}$, \eqref{TWT1} and \eqref{TWT2} hold.
Moreover, if the operator $\phi=B_1+B_2$ is linear, then $(\mathfrak{g},B_1,B_2)$ is a decomposable Rota-Baxter system of Lie algebras with $\phi$ idempotent.
\end{thm}

\proof
Recall that for $i=1,2$,  we have
$$\left.\frac{\mathrm{d}}{\mathrm{d}s }\right|_{s=0}\mathcal{B}_i(\exp(su))=\left.\frac{\mathrm{d}}{\mathrm{d}s }\right|_{s=0}\exp(sB_i(u))=B_i(u),\quad\forall u\in\mathfrak{g}.$$

We show the identity \eqref{TWT1}. In fact,
\begin{align*}
[B_1(u),B_1(v)]=&\left.\frac{\mathrm{d}^{2}}{\mathrm{d}s \mathrm{d}t}\right|_{s,t=0}\exp(sB_1(u))\exp(tB_1(v))\exp(-sB_1(u))\\
=&\left.\frac{\mathrm{d}^{2}}{\mathrm{d}s \mathrm{d}t}\right|_{s,t=0}\mathcal{B}_1(\exp(su))\mathcal{B}_1(\exp(tv))\mathcal{B}_1(\exp(su))^{-1}.
\end{align*}
Since $\mathcal{B}_1(\exp(su))^{-1}=\mathcal{B}_1(\mathcal{B}_1(\exp(su))^{-1}\mathcal{B}_2(\exp(su))^{-1})$,
\begin{align*}
 &[B_1(u),B_1(v)]\\
=&\left.\frac{\mathrm{d}^{2}}{\mathrm{d}s \mathrm{d}t}\right|_{s,t=0}\mathcal{B}_1(\exp(su))\mathcal{B}_1(\exp(tv))\mathcal{B}_1(\mathcal{B}_1(\exp(su))^{-1}\mathcal{B}_2(\exp(su))^{-1})\\
=&\left.\frac{\mathrm{d}^{2}}{\mathrm{d}s \mathrm{d}t}\right|_{s,t=0}\mathcal{B}_1(\exp(su))\mathcal{B}_1(\mathcal{B}_1(\exp(tv))\mathcal{B}_1(\exp(su))^{-1}\mathcal{B}_2(\exp(su))^{-1}\mathcal{B}_2(\exp(tv)))\\
=&\left.\frac{\mathrm{d}^{2}}{\mathrm{d}s \mathrm{d}t}\right|_{s,t=0}\mathcal{B}_1(\mathcal{B}_1(\exp(su))\mathcal{B}_1(\exp(tv))\mathcal{B}_1(\exp(su))^{-1}\\
&\qquad\qquad\qquad\qquad \times\mathcal{B}_2(\exp(su))^{-1}\mathcal{B}_2(\exp(tv))\mathcal{B}_2(\exp(su))).
\end{align*}
Using the Leibniz rule, we have
\begin{align*}
 [B_1(u),B_1(v)]=&B_1\left(\left.\frac{\mathrm{d}^{2}}{\mathrm{d}s \mathrm{d}t}\right|_{s,t=0}\mathcal{B}_1(\exp(su))\mathcal{B}_1(\exp(tv))\mathcal{B}_1(\exp(su))^{-1}\right.\\
+&\left.\left.\frac{\mathrm{d}^{2}}{\mathrm{d}s \mathrm{d}t}\right|_{s,t=0}\mathcal{B}_2(\exp(su))^{-1}\mathcal{B}_2(\exp(tv))\mathcal{B}_2(\exp(su)\right)\\
=&B_1([B_1(u),B_1(v)]-[B_2(u),B_2(v)]).
\end{align*}
Therefore we prove \eqref{TWT1}. The proof of \eqref{TWT2} is similar.

For any $t\in \mathbb{K}$ and $u\in\mathfrak{g}$, we have
\begin{equation*}
B_1(tu)=\left.\frac{\mathrm{d}}{\mathrm{d}s }\right|_{s=0}\mathcal{B}_1(\exp(s(tu)))=tB_1(u).
\end{equation*}
Therefore, to prove $B_1$ is linear, it's enough to show that $B_1(u)+B_1(v)=B_1(u+v)$ for any $u,v\in \mathfrak{g}$. Using Lemma \ref{LEM}, we have
\begin{equation*}
\begin{aligned}
B_1(u)+B_1(v)=&\left.\frac{\mathrm{d}}{\mathrm{d}s}\right|_{s=0}\mathcal{B}_1(\exp(su))\mathcal{B}_1(\exp(sv))\\
=&\left.\frac{\mathrm{d}}{\mathrm{d}s}\right|_{s=0}\mathcal{B}_1(\exp(su))\mathcal{B}_1(\Phi(\exp(sv))),
\end{aligned}
\end{equation*}
where $\Phi$ is the cocycle of the Rota-Baxter system $(G,\mathcal{B}_1,\mathcal{B}_2)$. Using Corollary \ref{B12}, we have $\Phi(1_G)=1_G$. Then
\begin{equation*}
\begin{aligned}
B_1(u)+B_1(v)=&\left.\frac{\mathrm{d}}{\mathrm{d}s}\right|_{s=0}\mathcal{B}_1(\mathcal{B}_1(\exp(su))\Phi(\exp(sv))\mathcal{B}_2(\exp(su)))\\
=&B_1(\phi(u)+\phi(v))=B_1(\phi(u+v)).
\end{aligned}
\end{equation*}
Note that using Lemma \ref{LEM}, we have
$$B_1(\phi(u))=B_1(u),\quad B_2(\phi(u))=B_2(u),\quad \forall u\in\mathfrak{g}.$$
Hence
$$B_1(u)+B_1(v)=B_1(u+v),$$
which implies that $B_1$ is linear. Similarly, we can show $B_2$ is linear. Hence $(\mathfrak{g},B_1,B_2)$ is a Rota-Baxter system of Lie algebras.

 Finally, by Lemma \ref{LIN} we see $(\mathfrak{g},B_1,B_2)$ is decomposable and $\phi$ is idempotent.
\qed

\begin{rem}
Without the condition that $\phi$ is linear, we can't derive that $B_1$ and $B_2$ are linear, which is different from the case of Rota-Baxter groups. That's because we can't ensure that $B_1(u)+B_1(v)=B_1(u+v)$ and $B_2(u)+B_2(v)=B_2(u+v)$ for any $u,v\in\mathfrak{g}.$ However, we have $B_1(u)+B_1(v)=B_1(\phi(u)+\phi(v))$ and $B_2(u)+B_2(v)=B_2(\phi(u)+\phi(v))$. Combining with Lemma \ref{LEM}, we obtain
$$B_1(\phi(u))+B_1(\phi(v))=B_1(\phi(u)+\phi(v)),\quad B_2(\phi(u))+B_2(\phi(v))=B_2(\phi(u)+\phi(v)).$$
Hence $B_1$ and $B_2$ are linear on $\phi(\mathfrak{g})$.
\end{rem}

\subsection{Twisted modified Yang-Baxter equations and Rota-Baxter systems of Lie algebras}

In this subsection, the notion of twisted modified Yang-Baxter equations is proposed. It's a generalization of the modified Yang-Baxter equation introduced in \cite{Semenov}. Its connection with the Rota-Baxter systems of Lie algebras is studied.

\begin{Def}
Let $\mathfrak{g}$ be a Lie algebra with linear operators $R:\mathfrak{g}\to \mathfrak{g}$ and $\phi:\mathfrak{g}\to \mathfrak{g}.$  Then pair $(R,\phi)$ is said to satisfy the twisted modified Yang-Baxter equations if
\begin{subequations}
\begin{align}
&[R(u),R(v)]+[\phi(u),\phi(v)]=R([R(u),\phi(v)]+[\phi(u),R(v)]),\label{TWT3}\\
&[R(u),\phi(v)]+[\phi(u),R(v)]=\phi([R(u),\phi(v)]+[\phi(u),R(v)])\label{TWT4}
\end{align}
\end{subequations}
hold for any $u,v\in \mathfrak{g}.$
\end{Def}

Note that if we take $\phi=\operatorname{id}$, then \eqref{TWT4} is always true and \eqref{TWT3} is just the modified Yang-Baxter equation \eqref{MYB} defined in \cite{Semenov}.

Recall that for a Lie algebra $\mathfrak{g}$, if the linear operator $R:\mathfrak{g}\to \mathfrak{g}$ satisfies the modified Yang-Baxter equation, then $P_1=\frac{1}{2}(R+\operatorname{id})$ is a Rota-Baxter operator of weight $-1$ and $P_2=\frac{1}{2}(R-\operatorname{id})$ is a Rota-Baxter operator of weight $1$. In the case for the modified Yang-Baxter equations, we have the following proposition.

\begin{prop}
Let $\mathfrak{g}$ be a Lie algebra with two linear operators $B_1:\mathfrak{g}\to\mathfrak{g}$ and $B_2:\mathfrak{g}\to\mathfrak{g}$. Then the following statements are equivalent:
\begin{description}
\item[(i)] $(\mathfrak{g},B_1,B_2)$ is a Rota-Baxter system of Lie algebras;
\item[(ii)] $(R,\phi)$ satisfies the twisted modified Yang-Baxter equations, where $\phi=B_1+B_2$ and $R=B_1-B_2$.
\end{description}
\end{prop}

\proof
For any $u,v\in\mathfrak{g}$, if $\phi=B_1+B_2$ and $R=B_1-B_2$, then
\begin{align*}
&[R(u),R(v)]+[\phi(u),\phi(v)]=2([B_1(u),B_1(v)]+[B_2(u),B_2(v)]),\\
&[R(u),\phi(v)]+[\phi(u),R(v)]=2([B_1(u),B_1(v)]-[B_2(u),B_2(v)]).
\end{align*}

Suppose $(\mathfrak{g},B_1,B_2)$ is a Rota-Baxter system of Lie algebras. Then using \eqref{TWT1} and \eqref{TWT2}, we have
\begin{align*}
[B_1(u),B_1(v)]+[B_2(u),B_2(v)]=&(B_1-B_2)([B_1(u),B_1(v)]-[B_2(u),B_2(v)])\\
=&\frac{1}{2}R([R(u),\phi(v)]+[\phi(u),R(v)])
\end{align*}
and
\begin{align*}
[B_1(u),B_1(v)]-[B_2(u),B_2(v)]=&(B_1+B_2)([B_1(u),B_1(v)]-[B_2(u),B_2(v)])\\
=&\frac{1}{2}\phi([R(u),\phi(v)]+[\phi(u),R(v)]).
\end{align*}
Now it is easy to get \eqref{TWT3} and \eqref{TWT4}.

Conversely, suppose that \eqref{TWT3} and \eqref{TWT4} hold. Then
\begin{align*}
&[B_1(u),B_1(v)]+[B_2(u),B_2(v)]=R([B_1(u),B_1(v)]-[B_2(u),B_2(v)]),\\
&[B_1(u),B_1(v)]-[B_2(u),B_2(v)]=\phi([B_1(u),B_1(v)]-[B_2(u),B_2(v)]).
\end{align*}
Hence we have
$$[B_1(u),B_1(v)]=\frac{R+\phi}{2}([B_1(u),B_1(v)]-[B_2(u),B_2(v)]),$$
which is \eqref{TWT1}, and
$$[B_2(u),B_2(v)]=\frac{R-\phi}{2}([B_1(u),B_1(v)]-[B_2(u),B_2(v)]),$$
which is \eqref{TWT2}.
\qed

\end{document}